\documentclass[12pt]{article}
\usepackage{amsmath,amsfonts}
\usepackage{mathrsfs}
\usepackage{color}
\usepackage{latexsym} 
\newcommand{\nref}[1]{(\ref{#1})}
\def\Cal{\cal}

\newtheorem{theorem}{Theorem}[section]
\newtheorem{remark}{Remark}[section]
\def\text#1{\hbox{#1}}

\def\endproof{\mbox{\ $\Box$}}
\def\Cal{\cal}

\newcommand{\N}        {{{\rm I \hspace{-.15em} N}}}
\newcommand{\R}        {{{\rm I\! R}}}

\def\1{\mbox{1\hspace{-.20em}I}}
\newcommand{\CZ}{{\Cal{Z}}}

\newcommand{\CX}{{\Cal{X}}}
\newcommand{\CA}{{\Cal{A}}}

\newcommand{\CN}{{\Cal{N}}}

\newcommand{\eq}{{ \,\stackrel{\Delta}{=}\, }}
\def\a{\alpha}
\def\t{\theta}

\def\e{\varepsilon}
\def\b{\beta}
\def\g{\gamma}
\def\la{\lambda}

\def\l{\left}
\def\r{\right}

\def\Var{\rm {Var}}

\def\CN{{\cal{N}}}

\def\CL{{\cal{L}}}

\def\Var{{\rm Var}}
\def\Cov{{\rm Cov}}

\newtheorem{proposition}[theorem]{Proposition}
\newtheorem{lemma}{Lemma}[section]
\newtheorem{defi}{Definition}[section]

\newtheorem{corollary}{Corollary}[section]
\voffset=-20mm \numberwithin{equation}{section}

\begin{document}
\title{\bf Detection boundary in sparse regression}
\author{
Ingster, Yu. I.
\thanks{St.Petersburg State Electrotechnical University, 5,
Prof. Popov str., 197376 St.Petersburg, Russia}, Tsybakov, A.B.
 \thanks{CREST, Timbre J340 3, av.\ Pierre Larousse, 92240
Malakoff cedex, France }, and Verzelen, N. {\thanks{INRA, UMR 729 MISTEA, 2
place Pierre Viala, F-34060 Montpellier, France}} }

\maketitle

\begin{abstract} We study
the problem of detection of a $p$-dimensional sparse vector
of parameters in the linear regression model with Gaussian noise.
We establish the detection boundary, i.e., the necessary and sufficient
conditions for the possibility of successful detection as both the
sample size $n$ and the dimension $p$ tend to the infinity. Testing procedures
that achieve this boundary are also exhibited. Our results encompass the
high-dimensional setting ($p\gg n$). The main
message is that, under some conditions, the detection boundary phenomenon
that has been proved for the Gaussian sequence model, extends to
high-dimensional linear regression. Finally, we establish the detection
boundaries when the variance of the noise is unknown. Interestingly, the
detection boundaries sometimes depend on the knowledge of the variance in a
high-dimensional setting.
\end{abstract}

\medskip





\section{Introduction}\label{SM}

 We consider the linear regression model with random
design:
\begin{equation}\label{Mod}
 Y_i=\sum_{j=1}^p\theta_j X_{ij}+\xi_i,\ i=1,...,n,
\end{equation}
where $\theta_j\in \R$ are unknown coefficients, $\xi_i$ are i.i.d.
$\CN(0,\sigma^2)$ random variables, $X_{ij}$ are random variables,
which are identically distributed, and $(X_{ij}, 1\le i\le n)$ are
independent for any fixed $j$ with $EX_{ij}=0,\ EX_{ij}^2=1$. We
study separately the settings with known $\sigma>0$ (then assuming
that $\sigma=1$ without loss of generality) and unknown $\sigma>0$.
We also assume that $X_{ij}$, $1\le j\le p, \ 1\le i\le n$, are
independent of $\xi_i, 1\le i\le n$.

Based on the observations $Z=(X,Y)$ where $X=(X_{ij}, 1\le j\le p, \
1\le i\le n)$, and $Y=(Y_i, 1\le i\le n)$, we consider the problem
of detecting whether the signal $\theta=(\theta_1,\dots,\theta_p)$
is zero (i.e., we observe the pure noise) or $\theta$ is some sparse
signal, which is sufficiently well separated from 0. Specifically,
we state this as a problem of testing the hypothesis $H_0:\theta=0$
against the alternative
$$
 H_{k,r}: \theta\in\Theta_k(r)=\{\theta\in \R^p_k: \|\theta\|\ge r\},
$$
where $\R^p_k$ denotes the $\ell_0$ ball in $\R^p$ of radius $k$,
$\|\cdot\|$ is the Euclidean norm, and $r>0$ is a separation
constant.

The smaller is $r$, the harder is to detect the signal. The question
that we study here is: What is the {\it detection boundary}, i.e.,
what is the smallest separation constant $r$ such that successful
detection is still possible? The problem is formalized in an
asymptotic minimax sense, cf. Section \ref{sec_problem} below. This
question is closely related to the previous work by several authors
on detection and classification boundaries for the Gaussian sequence
model
\cite{kjl07,dj04,dj08,dj09,hj10,Castro08,Castro10,I97,IS01a,IS01b,IS02a,IS02b,
IPT_arxiv,IPT_RS,JW,J03,J04}. These papers considered model
(\ref{Mod}) with $p=n$ and $X_{ij}=\delta_{ij}$, where $\delta_{ij}$
is the Kronecker delta, or replications of such a model (in
classification setting). The main message of the present work is
that, under some conditions, the detection boundary phenomenon
similar to the one discussed in those papers, extends to linear
regression. Our results cover the high-dimensional $p\gg n$ setting.

We now give a brief summary of our findings under the simplifying
assumption that all the regressors $X_{ij}$ are i.i.d. standard
Gaussian. We consider the asymptotic setting where
 $p\to\infty$, $n\to\infty$ and $k=p^{1-\beta}$ for some $\beta\in(0,1)$.
The results are different for moderately sparse alternatives
($0<\beta<1/2$) and highly sparse alternatives ($1/2<\beta<1$). We
show that for moderately sparse alternatives the detection boundary
is of the order of magnitude
\begin{equation}\label{rate1}
r\asymp \frac{p^{1/4}}{\sqrt{n}}\wedge \frac1{n^{1/4}},
\end{equation}
whereas for highly sparse alternatives ($1/2<\beta<1$) it is of the
order
\begin{equation}\label{rate2}
r\asymp \sqrt{\frac{k\log p}{n}}\wedge \frac1{n^{1/4}}.
\end{equation}
This solves the problem of optimal rate in detection boundary for
all the range of values $(p,n)$. Furthermore, for highly sparse
alternatives under the additional assumption
\begin{equation}\label{condc}
p^{1-\b}\log(p)=o(\sqrt{n})
\end{equation}
we obtain the sharp detection boundary,
i.e., not only the rate but also the exact constant. This sharp
boundary has the form
\begin{equation}\label{sep} r = \varphi(\beta)\sqrt{\frac{k\log
p}{n}},
\end{equation}
 where
\begin{equation}\label{def_phi}
     \varphi(\beta)=\begin{cases}\sqrt{2\b-1},& 1/2<\b\le 3/4,\\
     \sqrt{2}(1-\sqrt{1-\b}),& 3/4<\b<1.
     \end{cases}
\end{equation}
The function $\varphi(\cdot)$ here is the same as in the above
mentioned detection and classification problems, as first introduced
in \cite{I97}. We also provide optimal testing procedures. In
particular, the sharp boundary (\ref{sep})-(\ref{def_phi}) is
attained on the Higher Criticism statistic.

One of the applications of this result is related to transmission of
signals under compressed sensing, cf. \cite{d06,CanTao07}. Assume
that a sparse high-dimensional signal $\theta$ is coded using
compressed sensing with i.i.d. Gaussian $X_{ij}$ and then
transmitted through a noisy channel. Observing the noisy outputs
$Y_i$, we would like to detect whether the signal was indeed
transmitted. For example, this is of interest if several signals
appear in consecutive time slots but some slots contain no signal.
Then the aim is to detect informative slots. Our detection boundary
(\ref{sep}) specifies the minimal energy of the signal sufficient
for detectability. We note that $\varphi(\cdot)<\sqrt{2}$, so that
successful detection is possible for rather weak signals whose
energy is under the threshold $\sqrt{2k\log(p)/n}$. This can be
compared with the asymptotically optimal recovery of sparsity
pattern (RSP) by the Lasso in the same Gaussian model as ours
\cite{wainw09a,wainw09b}. Observe that the RSP property is stronger
than detection (i.e., it implies correct detection) but
\cite{wainw09a} defines the alternative by $\{\theta\in\R^p_k:
\,|\theta_j|\ge c\sqrt{\log(p)/n}, \forall j\}$ for some constant
$c>2$, which is better separated from the null than our alternative
$\Theta_k(r)$. Thresholds that are larger in order of magnitude are
required if one would like to perform detection based on estimation
of the values of coefficients in the $\ell_2$ norm
\cite{BRT09,CanTao07}.

In many applications, the variance of the noise $\xi$ is unknown.
Does the problem of detection become more difficult in this case? In
order to answer this question, we investigate the detection
boundaries in the unknown variance setting. Related work \cite{VV,V}
develop minimax bounds for detection in model (\ref{Mod}) under
assumptions different from ours and under unknown variance. However,
\cite{VV} does not provide a sharp boundary. Here, we prove that for
$\beta\in (1/2,1)$ and $k\log(p)\ll \sqrt{n}$, the detection
boundaries are the same for known and unknown variance. In contrast,
when $k\log(p)\gg \sqrt{n}$, the detection boundary is much larger
in the case of unknown variance than for known variance. We also
provide an optimal testing procedure for unknown variance.

After we have obtained our results, we became aware of the
interesting parallel unpublished work of Arias-Castro et
al.~\cite{arias}. There the authors derive the detection boundary in
model (\ref{Mod}) with known variance of the noise for both fixed
and random design. Their approach based on the analysis of the
Higher Criticism shares some similarities with our work. When the
variables $X_{ij}$ are i.i.d. standard normal and the variance is
known, we can directly compare our results with~\cite{arias}. In
\cite{arias} the detection boundaries analogous to (\ref{rate1}) and
(\ref{rate2}) do not contain the minimum with the $n^{-1/4}$ term,
because they are proved in a smaller range of values $(p,n)$ where
this term disappears. In particular, the conditions in~\cite{arias}
exclude the high-dimensional case $p\gg n$. We also note that, due
to the constraints on the classes of matrices $X$, \cite{arias}
obtains the sharp boundary (\ref{sep})-(\ref{def_phi}) under the
condition $p^{1-\b}(\log(p))^2=o(\sqrt{n})$ which is more
restrictive than our condition (\ref{condc}). The other difference
is that~\cite{arias} does not treat the case of unknown variance of
the noise.

\bigskip

Below we will use the following notation. We write $Z=(X,Y)$ where
$X=(X_{ij}, 1\le j\le p, \ 1\le i\le n)$, and $Y=(Y_i, 1\le i\le n)$
are the observations satisfying \nref{Mod}. Let $P_\theta$ be the
probability measure that corresponds to observations $Z$,
$P_{\theta,i}$ be those corresponding to observations
$Z_i=(X_{i1},...,X_{ip}, Y_i)$ with fixed $i=1,...,n$, and $P_X,
P_{X,i}$ be the probability measures corresponding to observations
$X$ or $X^{(i)}=(X_{i,j}, 1\le j\le p)$. We denote by
$P_{\theta}^{X}$ and $P_{\theta,i}^X$ the conditional distributions
of $Y$ given $X$ and of $Y_i$ given $X^{(i)}$ respectively. The
corresponding expectations are denoted by $E_{\theta}^{X}$ and
$E_{\theta,i}^X$. Clearly,
\begin{equation}\label{not1}
P_\theta(dZ)=P_{\theta}^{X}(dY) P_X(dX),\
P_{\theta,i}(dZ_i)=P_{\theta,i}^X(dY) P_{X,i}(dX^{(i)}),
\end{equation}
and
$$
 P_\theta(dZ)=\prod_{i=1}^nP_{\theta,i}(dZ_i).
$$
We denote by $X_j\in \R^n$ the $j$th   column of matrix
$X=(X_{ij})$, and set
$$
(X_{j},X_{l})=\sum_{i=1}^n X_{ij}X_{il},\quad \|X_j\|^2=(X_j,X_j).
$$

\section{Detection problem}\label{sec_problem}

For $\theta\in\R^p$, we denote by  $M(\theta)=\sum_{j=1}^p \1_{\{
\theta_j\not=0\}}$ the number of non-zero components of $\theta$,
where $\1_{\{A\}}$ is the indicator function. As above, let
$\R^p_k,\ 1\le k\le p,$ denote the $\ell_0$ ball in $\R^p$ of radius
$k$, i.e., the subset of $\R^p$ that consists of vectors $\theta$
with $M(\theta)\le k$, or equivalently, $\theta\in \R^p_k$ contains
no more than $k$ nonzero coordinates. In particular $\R^p_p=\R^p$.
Recall the notation $\Theta_k(r)=\{\theta\in \R^p_k: \|\theta\|\ge
r\}$.

We consider the problem of testing the hypothesis $H_0:\theta=0$
against the alternative $ H_{k,r}: \theta\in\Theta_k(r)$.  In this
paper we study the asymptotic setting where
 $p\to\infty$, $n\to\infty$ and $k=p^{1-\beta}$. The coefficient  $\beta\in
[0,1]$ is called the sparsity index. 
We assume in this section that $\sigma^2$ is known. Modifications
for the case of unknown variance are discussed in Section
\ref{subsec_problem_statement_unknown}. 

We call a test any measurable function $\psi(Z)$ with values in
$[0,1]$. For a test $\psi$, let
$
\a(\psi)=E_0(\psi)
$
be the type I error probability and
$
\b(\psi,\t)=E_\t(1-\psi)
$
be the type II error probability for the alternative
$\t\in\Theta\subset \R^p$. We set
$$
\b(\psi)=\b(\psi,\Theta)=\sup_{\t\in\Theta}\b(\psi,\t),\quad
\g(\psi)=\g(\psi,\Theta)=\a(\psi)+\b(\psi,\Theta) \ .
$$
We denote by
$\b(\a)=\b_{n,p,k}(\a,r)$ the minimax type II error probability for
a given level $\a\in (0,1)$,
$$
\b(\a)=\inf_{\psi:\a(\psi)\le\a}\b(\psi,\Theta_k(r)),\quad 0\le
\b(\a)\le 1-\a\ .
$$
Accordingly, we denote by $\gamma=\gamma_{n,p,k}(r)$ the minimax
total error probability in the hypothesis testing problem:
$$
\gamma=\inf_{\psi}\g(\psi,\Theta_k(r)),
$$
where the infimum is taken over all tests $\psi$. Clearly,
$$
\gamma=\inf_{\a\in (0,1)}(\a+\b(\a)),\quad 0\le \g\le 1.
$$

The aim of this paper is to establish the asymptotic detection
boundary, i.e., the conditions on the separation constant
$r=r_{n,p,k}$, which delimit the zone where $\gamma_{n,p,k}(r)\to 1$
(indistinguishability) from that where $\gamma_{n,p,k}(r)\to 0$
(distinguishability). The distinguishability is equivalent to
$\b(\a)\to 0,\ \forall \ \a\in (0,1)$. We are interested in tests
$\psi=\psi_{n,p}$ or $\psi_\a=\psi_{n,p,\a}$ such that either $
\g(\psi)\to 0$ or $\a(\psi_\a)\le \a+o(1)$, and $\b(\psi_\a)\to 0$.
Here and later the limits are taken as $p\to\infty,\ n\to\infty$
unless otherwise stated.

\section{Assumptions on $X$}\label{sec_assump}

We will use at different instances some of the following conditions
on the random variables $X_{ij}$.

\smallskip

{\it

\noindent{\bf A1}. The random variables $X_{ij}$ are uncorrelated,
i.e., $EX_{ij}X_{il}=0$ for all $1\le j<l\le p$.

\noindent{\bf A2}. The random variables $X_{ij}$, $1\le j\le p, \
1\le i\le n$, are independent.

\noindent{\bf A3}. The random variables $X_{ij}$, $1\le j\le p, \
1\le i\le n$, are i.i.d. standard Gaussian: $X_{ij}\sim \CN(0,1)$.

}

\smallskip

Let $U_j,\ 1\le j\le p$ be random variables such that we have
equality in distribution $\CL(U_j,U_l)=\CL(X_{ij},X_{il}),\quad 1\le
j<l\le p$. We will need the following technical assumptions.

{\it

\begin{eqnarray}\label{moments}
{\bf B1}.\hspace{4.2cm}\max_{1\le j<l\le p}E((U_jU_l)^4)=O(1)\ . \hspace{4.2cm}
\end{eqnarray}

\noindent
{\bf B2}. There exists $h_0>0$ such that $\max_{1\le j\le l\le
p}E(\exp(hU_jU_l))<\infty$ for $|h|<h_0$, and
\begin{equation}\label{TA1}
\log^3(p)=o(n).
\end{equation}

\noindent
{\bf B3}. There exists $m\in \N$ such that $\max_{1\le j\le l\le
p}E(|U_jU_l|^m)<\infty$, and
\begin{equation}\label{TA2}
\log^2(p)p^{4/m}=o(n).
\end{equation}

}

\noindent Assumption {\bf B1} implies that
\begin{eqnarray}\label{moments0}
\max_{1\le j<l\le p}E(|U_jU_l|^m)&=&O(1),\quad m=2,3,4\ .
\end{eqnarray}
In particular, Assumption {\bf B1} holds true under {\bf A2} if
\begin{equation}\label{moments2}
\max_{1\le j\le p}E(U_j^4)=O(1).
\end{equation}
If $(X_{ij}X_{il}, i=1,\dots,n)$ are independent zero-mean random
variables, we have (cf. \cite{petrov}, p. 79):
$$
E|(X_{j},X_{l})|^m\le C(m)n^{m/2-1}\sum_{i=1}^n
E(|X_{ij}X_{il}|^m),\quad m>2\ .
$$
This and \nref{moments0} yield
\begin{eqnarray}
\sum_{1\le j<l\le p}E(|(X_j,X_l)|^m)&=&O(n^{m/2}p^2),\quad m=2,3,4\ .
\label{moments1}
\end{eqnarray}

\noindent
Finally,  Assumptions {\bf B1} and {\bf B2} hold true under {\bf A3}
and \nref{TA1}.

\section{Main results}\label{D}

\subsection{Detection boundary under known
variance}\label{section_known} For this case  we suppose $\sigma=1$
without loss of generality.

\subsubsection{Lower
bounds}\label{DL} We first present the lower bounds on the detection
error, i.e., the indistinguishability conditions. We assume that
$k=p^{1-\b}, \b\in (0,1)$. Indistinguishability conditions consist
of two joint conditions on the radius $r=r_{np}$.  The first one is
\begin{equation}\label{L1}
r_{np}^2=o(n^{-1/2}).
\end{equation}
The second condition differs according to whether $\b\le 1/2$ or
$\b> 1/2$. If $\beta\le 1/2$ (i.e. $p=O(k^2)$),  which corresponds to moderate
sparsity, we require that
\begin{equation}\label{L2a}
r_{np}^2=o(\sqrt{p}/n)\ .
\end{equation}
The case $\beta>1/2$ (i.e. $k^2=o(p)$) corresponds to high
sparsity. We define $x_{np}$ by $r_{np}=x_{n,p}\sqrt{k\log(p)/n}$.
Then, we require that
\begin{equation}\label{L2b}
           \lim\sup (x_{n,p}-\varphi(\beta))<0,
\end{equation}
     where $\varphi(\beta)$ is defined in (\ref{def_phi}).
     Clearly, condition \nref{L2b} implies $r_{np}^2=O(k\log(p)/n)$,
which is stronger than \nref{L2a} when $\beta>1/2$.

\begin{theorem}\label{TL} Assume {\bf A1}, {\bf B1}, $k=p^{1-\b}$
and  either  {\bf B2} or {\bf B3}. We also require that $r_{np}$ satisfies
\nref{L1} and
either
\nref{L2a} (for $\b\in (0,1/2]$) or \nref{L2b} (for $\b\in
(1/2,1)$). Then, asymptotic distinguishability is impossible, i.e.,
$\g_{n,p,k}(r_{np})\to 1$.
\end{theorem}

\begin{remark}\label{RL1}
{This theorem can be extended to non-random design matrix $X$.
Inspection of the proof shows that, instead of {\bf B1}, we only
need the assumption: For some $B_{n,p}$ tending to $\infty$ slowly
enough,
\begin{equation}\label{b11}
\sum_{1\le j<l\le p}|(X_j,X_l)|^m<B_{n,p}n^{m/2}p^2,\quad m=2,3,4.
\end{equation}
Indeed, {\bf B1} is used in the proofs only to assure that
\nref{b11} holds true with $P_X$-probability tending to 1 (this is
deduced from assumption {\bf B1} and \nref{moments1}).

Also instead of {\bf B2} and {\bf B3}, we can assume that there
exists $\eta_{n,p}\to 0$ such that
\begin{equation}\label{b22}
r_{np}^2\max_{1\le j<l\le p}|(X_j,X_l)|<\eta_{n,p}k,\quad \max_{1\le
j\le p}|\|X_j\|^2-n|<\eta_{n,p}n.
\end{equation}
Under {\bf B2}, {\bf B3}, relations \nref{b22} hold with
$P_X$-probability tending to 1, see  Corollary \ref{C1}.

The result of the theorem remains valid for non-random matrices $X$
satisfying \nref{b11} and \nref{b22}.
 }
\end{remark}

 \subsubsection{Upper bounds}\label{U}

 In order to construct a test procedure that achieves the detection
boundary,
 we combine several tests.\\

\noindent First, we study the widest non-sparse case $k=p$, i.e., we
consider $\Theta_p(r)=\{\theta\in \R^p: \|\theta\|\ge r\}$. Consider
the statistic
\begin{equation}\label{stat0}
                  t_0=(2n)^{-1/2}\sum_{i=1}^n(Y_i^2-1),
\end{equation}
which is the $H_0$-centered and normalized version of the classical
$\chi^2_n$-statistic $\sum_{i=1}^nY_i^2$. The corresponding tests
$\psi_\a^0$ and $\psi^0$ are of the form:
$$
\psi_\a^0=\1_{t_0>u_\a},\quad \psi^0=\1_{t_0>T_{np}}
$$
where $\a\in (0,1)$, $u_\a$ is the $(1-\a)$-quantile of the standard
Gaussian distribution and $T_{np}$ is any sequence satisfying
$T_{np}\to\infty$.

\begin{theorem}\label{TU1}For all $\a\in
(0,1)$, we have:

(i) Type I errors satisfy $ \a(\psi_\a^0)=\a+o(1)$ and $\a(\psi^0)=o(1)$.

(ii) Type II errors. Assume {\bf A2} and {\bf B1}  and consider a radius
$r_{np}$ such that
~\\$nr^4_{np}\to\infty$. Then, we have $\b(\psi_\a^0, \Theta_p(r_{np}))\to 0$.
If
$T_{np}$ is chosen such that $\lim\sup
T_{np}n^{-1/2}r_{np}^{-2}<1 $, then
$\b(\psi^0, \Theta_p(r_{np}))\to 0$.
\end{theorem}
Recall that we can replace {\bf B1} by \nref{moments2} under {\bf
A2}. If $nr^4_{np}\to \infty$, then
one can take $T_{np}$ such that
$\g_{n,p}(\psi^0,\Theta_p(r_{np}))\to 0$ under {\bf A2}, {\bf B1}. This upper
bound corresponds to the part \nref{L1} of the detection boundary.\\

%


We now introduce a test $\psi_{\alpha}^1$ that achieves the second
boundary (\ref{L2a}). Consider the following kernel
$$
K(Z_i,Z_k)=p^{-1/2}Y_iY_k\sum_{j=1}^p X_{ij}X_{kj}. \
$$
The $U$-statistic $t_1$ based on the kernel $K$ is defined by
$$
t_1=N^{-1/2}\sum_{1\le i<k\le n}K(Z_i,Z_k),\quad N={n(n-1)}/{2}.
$$
Note that the  $U$-statistic $t_1$ can be viewed as the
$H_0$-centered and normalized version of the statistic
$\chi^2_p=n\sum_{j=1}^p\hat\theta_j^2$ based on the estimators
$\hat\theta_j=n^{-1}\sum_{i=1}^nY_iX_{ij}$:
$$
\chi^2_p=2n^{-1}\sum_{j=1}^p\sum_{1\le i<k\le
n}Y_iY_kX_{ij}X_{kj}+n^{-1}\sum_{j=1}^p\sum_{i=1}^nY_i^2X_{ij}^2.
$$
Indeed, up to a normalization, the first sum is the $U$-statistic
$t_1$, and moving off the second sum corresponds to centering.

 Given $\a\in (0,1)$, we consider the test $\psi_\a^1=\1_{t_1>u_\a}$.

\begin{theorem}\label{TU2} Assume {\bf A2} and
{\bf B1}. For all $\a\in (0,1)$ we
have:

(i) Type I error satisfies: $ \a(\psi_\a^1)=\a+o(1).$

(ii) Type II errors. Assume that $p=o(n^2)$ and
consider a radius $r_{np}$ such that
$nr^2_{np}/\sqrt{p}\to\infty$.
Then,
$
\b(\psi_\a^1, \Theta_p(r_{np}))\to 0
$.
\end{theorem}


\begin{remark}\label{RU3}
 Combining the tests $\psi^0_{\a/2}$ and $\psi^1_{\a/2}$
 we obtain the test $\psi^*_\a=\max(\psi^0_{\a/2},\psi^1_{\a/2})$ of asymptotic
 level not more than $\a$. Moreover, it achieves $\b(\psi_\a^1,
\Theta_p(r_{np}))\to 0$ for any radius $r_{np}$ satisfying
$$r^2_{np}\gg \frac{\sqrt{p}}{n}\wedge \frac{1}{\sqrt{n}}\ .$$
 We can omit the condition $p=o(n^2)$ since the test
$\psi^0_{\a/2}$ achieves the optimal rate for $p\geq n$. Combining
this bound with Theorem \ref{TL}, we conclude that $\psi^*_\a$
simultaneously achieves the optimal detection rate for all $\beta\in
(0,1/2]$.
\end{remark}

We now turn to testing in the highly-sparse case $\beta\in(1/2,1)$.
Here we use a version of "Higher Criticism Tests" (HC-tests, cf.
\cite{dj04}). Set
$$
y_i=(Y,X_i)/\|Y\|,\quad 1\le i\le p.
$$
Let $q_i=P(|\CN(0,1)|>|y_i|)$ be the $p$-value of the $i$-th
component and let $q_{(i)}$ denote these quantities sorted in
\emph{increasing order}. We define the HC-statistic by
\begin{eqnarray}\label{definition_HC}
 t_{HC}= \max_{i,\  q_{(i)}\leq
1/2}\frac{\sqrt{p}(i/p-q_{(i)})}{\sqrt{q_{(i)}(1-q_{(i)})}}\ .
\end{eqnarray}
Given a constant $a>0$, the HC-test $\psi^{HC}$ rejects $H_0$
when the
statistic $t_{HC}$ is larger than $(1+a)\sqrt{2\log\log p}$.

\begin{remark}
The cutoff $1/2$ in the definition \nref{definition_HC} of $t_{HC}$
can be replaced by any $c\in (0,1)$.
\end{remark}

\begin{theorem}\label{TU4} Assume {\bf A3} ($X_{ij}$ are i.i.d.
standard Gaussian).

(i) Type I error satisfies $ \a(\psi^{HC})=o(1).$

(ii) Type II error. Consider  $k=p^{1-\b}$ with $\b\in (1/2,1)$
and assume that
$k\log(p)=o(n)$. Take a radius $r_{np}=x_{np}\sqrt{k\log(p)/n}$ such that
$ \lim\inf(x_{np}-\varphi(\b))>0$. Then, we have  $
\b(\psi^{HC},\Theta_{k}(r_{np}))\to
0 $.

\end{theorem}

\begin{remark}
If $k\log(p)=o(n)$, the HC-test asymptotically detects any $k$-sparse signal
whose rescaled intensity $r_{np}\sqrt{n/(k\log(p))}$ is above the detection
boundary $\varphi(\beta)$.
\end{remark}

\begin{remark}
Assume ${\bf A3}$.
 Combining the tests $\psi^0_{\a}$ and $\psi^{HC}$, we obtain the test
$\psi^{*,HC}_\a=\max(\psi^0_{\a},\psi^{HC})$ of asymptotic level not more than
$\a$. Moreover, it achieves $\b(\psi_\a^1,
\Theta_k(r_{np}))\to 0$ for any radius $r_{n,p}$ satisfying $$\lim\inf
x_{n,p}\geq \varphi(\beta)\quad \text{ or }\quad r^2_{np}\gg
\frac{1}{\sqrt{n}}\ .$$
 We can omit the condition $k\log(p)=o(n)$ since the test
$\psi^0_{\a}$ achieves the optimal rate for $k\log(p)\gg \sqrt{n}$.
Combining this bound with Theorem \ref{TL}, we conclude that
$\psi^{*,HC}_{\a}$ simultaneously achieves the optimal detection
rate for all $\beta\in (1/2,1)$.

\end{remark}

In conclusion, under Assumption ${\bf A3}$, the test
$\max(\psi^0_{\a/2},\psi^1_{\a/2},\psi^{HC})$ simultaneously
achieves the optimal detection rate for all $\beta\in (0,1)$. The
detection boundary is of the order of magnitude
\begin{equation}\label{om}
r\asymp \sqrt{\frac{k\log p}{n}}\wedge \frac{1}{n^{1/4}}\ .
\end{equation}
Furthermore, we establish the sharp detection boundary (i.e., with
exact asymptotic constant) of the form
$$r=\varphi(\beta)\sqrt{\frac{k\log p}{n}}$$
for $\beta>1/2$ and $k\log(p) = p^{1-\beta}\log(p) = o(\sqrt{n})$.

\subsection{Detection boundary under unknown
variance}\label{section_unknown}

\subsubsection{Detection
problem}\label{subsec_problem_statement_unknown}
Since the
variance  of the noise  is now assumed to be unknown,  the tests
$\psi$ under study should not require the knowledge of
$\sigma^2$. The type I error  probability is now taken uniformly over
$\sigma>0$:
\begin{eqnarray*}
 \alpha^{un}(\psi) = \sup_{\sigma>0}E_{0,\sigma}(\psi)\ .
\end{eqnarray*}
The type II error probability over an alternative $\Theta\subset \mathbb{R}^p$
is
\begin{eqnarray}\label{definition_beta}
 \beta^{un}(\psi,\Theta) =
\sup_{\theta\in\Theta,\sigma>0}\beta(\psi,\theta\sigma,\sigma)=\sup_{
\theta\in\Theta,\sigma>0} E_{\theta\sigma,\sigma}(1-\psi)\ .
\end{eqnarray}
Similarly to the setting with known variance, we consider the sum of
the two errors: $$\gamma^{un}(\psi,\Theta)= \alpha^{un}(\psi)+
\beta^{un}(\psi,\theta).$$
 Finally, the minimax total error probability in the
hypothesis testing problem
 with unknown variance is
$$ \gamma_{n,p,k}^{un}(r) = \inf_{\psi} \gamma^{un}(\psi,\Theta_k(r))$$

\subsubsection{Lower bounds}\label{section_lower_unknown}

Take $r_{np}= x_{n,p}\sqrt{k\log(p)/n}$. As in the case of known
variance, we consider the condition
\begin{eqnarray}\label{condition_separation}
 \lim\sup (x_{n,p}-\varphi(\beta))<0 \ .
\end{eqnarray}

\begin{theorem}\label{prte_lower_bound_unknown}
Fix some $\beta>1/2$ and assume ${\bf A3}$.  If Condition
(\ref{condition_separation}) holds and if $k\log(p)= o(n)$, then
distinguishability is impossible, i.e., $
\gamma^{un}_{n,p,k}(r_{np})\to 1\ .$

If $k\log(p)/n\rightarrow \infty$, then for any radius $r>0$,
distinguishability is impossible, i.e.
$\gamma^{un}_{n,p,k}(r)\rightarrow\  1$.
\end{theorem}

The detection boundary stated in Theorem \ref{prte_lower_bound_unknown} does not
depend on the unknown $\sigma^2$. This is due to the definition
(\ref{definition_beta}) of the type II error
probability $\beta^{un}(\psi,\Theta_k(r))$ that considers alternatives of the
form $\sigma\theta$ with $\theta\in \Theta_k(r)$.

\subsubsection{Upper bounds}\label{section_upper_unknown}

The HC-test $\psi^{HC}$ defined in (\ref{definition_HC}) still achieves the
optimal detection rate when the variance is unknown as shown by the next
proposition.

\begin{proposition}\label{prte_upper_bound_unknown}
 Assume {\bf A3} ($X_{ij}$ are i.i.d.
standard Gaussian).

(i) Type I error satisfies $ \a^{un}(\psi^{HC})=o(1).$

(ii) Type II error. Consider  $k=p^{1-\b}$ with $\b\in (1/2,1)$
and assume that
$k\log(p)=o(n)$. Take a radius $r_{np}=x_{np}\sqrt{k\log(p)/n}$ such that
$ \lim\inf(x_{np}-\varphi(\b))>0$. Then, we have  $
\b^{un}(\psi^{HC},\Theta_{k}(r_{np}))\to
0 $.
\end{proposition}

In conclusion, in the setting with unknown variance we prove that
the sharp detection boundary (i.e., with exact asymptotic constant)
of the form
$$\varphi(\beta)\sqrt{\frac{k\log
p}{n}}$$ holds for $\beta>1/2$ and $k\log(p)=p^{1-\b}\log(p)=o(n)$,
i.e., for a larger zone of values $(p,n)$ than for the case of known
variance. However, this extension corresponds to $(p,n)$ for which
the rate itself is strictly slower than under the known variance.
Indeed, if the variance $\sigma^2$ is known, as shown in
Section~\ref{section_known}, the detection boundary is of the order
(\ref{om}). Thus, there is an asymptotic difference in the order of
magnitude of the two detection boundaries for $k\log(p)\gg
\sqrt{n}$.

\section{Proofs of the lower bounds}\label{PL}
\subsection{The prior}
Take $c\in (0,1),\ h=ck/p,\ b=r_{np}/c\sqrt{k},\ a=b\sqrt{n}$. Note
that the condition $r_{np}^2=o(1/\sqrt{n})$ is equivalent to
$b^4k^2n=o(1)$.

Let us consider a random vector $\theta=(\theta_j)$ with coordinates
$$
\theta_j=b\e_j,\quad\text{where}\quad \e_j\in (0,+1,-1)\quad
\text{iid},
$$
such that
$$
Prob(\e_j=0)=1-h,\quad Prob(\e_j=+1)=Prob(\e_j=-1)=h/2.
$$
This introduces a prior probability measure $\pi_j$ on $\theta_j$
and the product prior measure $\pi=\prod_{j=1}^p\pi_j$ on $\theta$.
The corresponding expectation and variance operators will be denoted
by $E_{\pi}$ and $\Var_{\pi}$.

\begin{lemma}\label{Lem1} Let $k\to\infty$. Then
$$\pi(\theta\in \R_k^p,\ \|\theta\|\ge r)\to 1.$$
\end{lemma}
{\bf Proof}. Observe that
$$
\|\theta\|^2=b^2\sum_{j=1}^p\e_j^2,\ m(\theta)=\sum_{j=1}^p|\e_j|.
$$
We have
$$
E_{\pi}(\|\theta\|^2)=b^2ph=r^2/c,\quad E_{\pi}(m(\theta))=ph=ck,
$$
and
$$
\Var_{\pi}(\|\theta\|^2)\le ph b^4=r^4/(kc^3),\quad
\Var_{\pi}(m(\theta))\le ph=ck.
$$
Applying the Chebyshev inequality, we get with $C=c^{-1}>1$,
$$
\pi(\|\theta\|^2<r^2)=\pi(E_{\pi}(\|\theta\|^2)-\|\theta\|^2>r^2(C-1))
\le\frac{\Var_{\pi}(\|\theta\|^2)}{r^4(C-1)^2}\to
0,
$$
and similarly,
$
\pi(m(\theta)>k)\to 0.
$
\endproof

Lemma \ref{Lem1} implies  that, in order to obtain asymptotic lower
bounds for the minimax problem, we only have to study the Bayesian
problem which corresponds to the prior $\pi$, see for instance
\cite{IS02a}, Proposition 2.9. Consider the mixture
$$
P_\pi(dZ)=E_\pi P_\theta(dZ)=\int_{\R^p}P_\theta(dZ)\pi(d\theta)
$$
and the likelihood ratio
$$
L_\pi(Z)=\frac{dP_\pi}{dP_0}(Z).
$$
In order to prove the lower bounds we only need to check that
\begin{equation}
L_\pi(Z)\to 1\quad\text{in}\ P_0-\text{probability}.
\end{equation}

Consider $x=\lim\sup x_{n,p}$. If $\beta\le 1/2$, then $x=0$ since
$nb^2=O(1)$. For $\beta>1/2$, we take $c\in (0,1)$ such that
$x_c=x/c<\varphi(\b)$, which is possible as $x<\varphi(\b)$.  We
will use the short notation $x$ and $a$ for $x_c$ and
$a_c=b\sqrt{n}=a/c$. We set
$$
a_j=b\|X_j\|,\quad y'_j=(X_j,Y)/\|X_j\|,\quad
x_j=a_j/\sqrt{\log(p)},\quad T_j=a_j/2+\log(h^{-1})/a_j,
$$
which corresponds to $he^{-a_j^2+a_jT}=1$.

\subsection{Study of the likelihood ratio $L_\pi$}

First observe that by \nref{not1}
$$
P_\pi(dZ)=P_X(dX)E_\pi \l(P_{\theta}^X(dY)\r),\quad
L_\pi(Z)=E_\pi\l(\frac{dP_{\theta}^X}{dP_{0}^X}(Y)\r).
$$
Note that conditional measure $P_{\theta}^X$ corresponds to
observation of the Gaussian vector $\CN(v,I_n)$ where
$v=\sum_{j=1}^p\theta_jX_j$, $I_n$ is the $n\times n$ identity
matrix, and the likelihood ratio under the expectation is
$$
\frac{dP_{\theta}^X}{dP_{0}^X}(Y)=\exp(-\|v\|^2/2+(v,Y))=
g_\theta(Z)e^{-\Delta(X,\theta)},
$$
where
\begin{equation}\label{gD}
g_\theta(Z)=\prod_{j=1}^p\exp(-\theta_j^2\|X_j\|^2/2+\theta_j(X_j,Y)),\quad
\Delta(X,\theta)=2\sum_{1\le j<l\le p}\theta_j\theta_l(X_j,X_l).
\end{equation}
Put
$$
\Lambda(Z)=E_\pi
\l(g_\theta(Z)\r)=\prod_{j=1}^p(1-h+he^{-b^2\|X_j\|^2/2}\cosh(b(X_j,Y)))\ .
$$
We define $\eta_j= e^{-b^2\|X_j\|^2/2}\cosh(b(X_j,Y))-1$.
Take now $\delta>0$ and introduce the set
$$
\Sigma_X=\{\theta\in\R^p:|\Delta(X,\theta)|\le \delta\}.
$$
We can write
$$
L_\pi(Z)=\int_{\R^p}g_\theta(Z)e^{-\Delta(X,\theta)}\pi(d\theta)\ge
e^{-\delta}\int_{\Sigma_X}g_\theta(Z)\pi(d\theta)=e^{-\delta}
\Lambda(Z)\pi_Z(\Sigma_X),
$$
where $\pi_Z=\prod_{j=1}^p\pi_{Z,j}$ is the random probability
measure  on $\R^p$ with the density
$$
\frac{d\pi_Z}{d\pi}(\theta)=\frac{g_\theta(Z)}{\Lambda(Z)}=
\prod_{j=1}^p\frac{d\pi_{Z,j}}{d\pi_j}(\theta);\quad
\frac{d\pi_{Z,j}}{d\pi_j}(\theta)=
\frac{e^{-\theta_j^2\|X_j\|^2/2+\theta_j(X_j,Y)}}{1+h\eta_j},\quad
\theta_j\in\{0,\pm b\},
$$
i.e., the measure $\pi_{Z,j}$ is supported at the points $\{0, b,
-b\}$ and
$$
\pi_{Z,j}(0)=\frac{1-h}{1+h\eta_j},\ \pi_{Z,j}(\pm
b)=\frac{h_{Z,j}^{\pm}}{2},\quad
h_{Z,j}^{\pm}=\frac{he^{{d_j}^{\pm}}}{1+h\eta_j},
$$
where we set
$$
d_j^{\pm}=-a_j^2/2\pm a_jy'_j\ ,\quad
\eta_j=\frac{e^{d^+_j}}{2}+\frac{e^{d^-_j}}{2}-1 .
$$

\begin{proposition}\label{P1} In $P_{0}$-probability,
\begin{equation}\label{probab*}
\pi_Z(\Sigma_X)\to 1.
\end{equation}
\end{proposition}
{\bf Proof of Proposition \ref{P1}} is given in Section
\ref{ProofP1}.

\begin{proposition}\label{P2} In $P_{0}$-probability,
\begin{equation}\label{prob2}
\Lambda(Z)\to 1.
\end{equation}
\end{proposition}
{\bf Proof of Proposition \ref{P2}} is given in Section
\ref{ProofP2}.\\

Propositions \ref{P1} and \ref{P2} imply that, for any $\delta>0$,
$$
P_0(Z: L_\pi(Z)>1-\delta)\to 1.
$$
Since $E_{0}L_\pi=1$ and $L_\pi(Z)\ge 0$, this yields $L_\pi\to 1$ in
$P_{0}$-probability. This yields indistinguishability in the problem.
\endproof

\subsection{Proof of Proposition \ref{P1}}\label{ProofP1}

\subsubsection{Replacing the measure $\pi_Z$ by
$\tilde\pi_Z$}\label{SS2}

Let us consider the random measure $ \tilde \pi_Z=\prod_{j=1}^p
\tilde \pi_{Z,j}, $ where $\tilde \pi_{Z,j}$ is supported at the
points $\{0, b, -b\}$ and
$$
\tilde\pi_{Z,j}(0)=1-\frac{q_{Z,j}^+}{2} -\frac{q_{Z,j}^-}{2},\
\tilde\pi_{Z,j}(\pm
b)=\frac{q_{Z,j}^{\pm}}{2},
$$
where
$$
q_{Z,j}^{\pm}=(h/2)e^{{d_j}^{\pm}}\1_{\CA_j^{\pm}},\quad
\CA_j^{\pm}=\{he^{{d_j}^{\pm}}<1\}=\{\pm y'_j<T_j\}
$$
and observe that the event $\CA_j^{\pm}$ implies $q_{Z,j}^{\pm}\le
1/2$, i.e, the measures $\tilde\pi_{Z,j}$ are correctly defined. We define the
event $\CA= \CA_{np}= \cap_{j=1}^p (\CA_j^+ \cap \CA_j^-)$.

\begin{lemma}\label{L0}
$$
P_0(\CA_{n,p})\to 1.
$$
\end{lemma}
{\bf Proof}. Denote $A^c$ the complement of the event $A$. Since
$y'_j\sim \CN(0,1)$ under $P_0$, we have
$$
P_0^X((\CA_{n,p})^c)\le \sum_{j=1}^p
P_0^X((\CA_j^+)^c)+P_0^X((\CA_j^-)^c) =2\sum_{j=1}^p\Phi(-T_j).
$$
By  Corollary \ref{C1} we get $a_j=b\|X_j\|\sim b\sqrt{n}$ uniformly in
$1\le j\le p$ in $P_X$-probability. By \nref{Phi0} this implies
$\sum_{j=1}^p\Phi(-T_j)=o(1)$  in $P_X$-probability.
\endproof

\medskip

We can replace the measure $\pi_Z$ by $\tilde\pi_Z$ in
\nref{probab*}. This follows from the following lemma

\begin{lemma}\label{Lem2} In $P_{0}$-probability,
\begin{equation}\label{prob1}
E_{\tilde\pi_Z}|d\pi_Z/d\tilde\pi_Z-1|\to 0.
\end{equation}
\end{lemma}
{\bf Proof}. Applying the equality
$E_{\tilde\pi_Z}(d\pi_Z/d\tilde\pi_Z)=1$ and the inequality $1+x\le e^x$, we get
\begin{eqnarray*}
(E_{\tilde\pi_Z}|d\pi_Z/d\tilde\pi_Z-1|)^2&\le&
E_{\tilde\pi_Z}(d\pi_Z/d\tilde\pi_Z-1)^2=E_{\tilde\pi_Z}
(d\pi_Z/d\tilde\pi_Z)^2-1\\
&=&
\prod_{j=1}^pE_{\tilde\pi_{Z,j}}\l(d\pi_{Z,j}/d\tilde\pi_{Z,j}\r)^2-1\\
&=&\prod_{j=1}^p(1+E_{\tilde\pi_{Z,j}}(d\pi_{Z,j}/d\tilde\pi_{Z,j}-1)^2)-1\\
&\le&\exp\l(\sum_{j=1}^pE_{\tilde\pi_{Z,j}}(d\pi_{Z,j}/d\tilde\pi_{Z,j}-1)^2\r)
-1.
\end{eqnarray*}
Consequently, we only have to prove that in $P_{0}$-probability,
$$
H(Z)=\sum_{j=1}^pE_{\tilde\pi_{Z,j}}(d\pi_{Z,j}/d\tilde\pi_{Z,j}-1)^2\to
0.
$$
Since $H(Z)\ge 0$, the last relation follows from
$$
E_0^X(H)\to 0,\quad \text{in $P_X$-probability}
$$
by Markov inequality. Observe that
$$
E_{\tilde\pi_{Z,j}}(d\pi_{Z,j}/d{\tilde\pi_{Z,j}}-1)^2=
\frac{(h_{Z,j}^+-q_{Z,j}^+)^2}{2q_{Z,j}^+}+
\frac{(h_{Z,j}^--q_{Z,j}^-)^2}{2q_{Z,j}^-}
+\frac{(h_{Z,j}^++h_{Z,j}^--q_{Z,j}^+-q_{Z,j}^-)^2}{2(2-q_{Z,j}^+-q_{Z,j}^-)}.
$$
By Lemma \ref{L0}, it is sufficient to study these terms under the event
$\CA$ which corresponds to $\max_{1\le j\le p}q_{Z,j}^{\pm}\le 1/2$.
Under this event, we have $h_{Z,j}^\pm=q_{Z,j}^\pm/\la_j,\
\la_j=1+q_{Z,j}^++q_{Z,j}^--h$, and direct calculation gives
\begin{eqnarray*}
\frac{(h_{Z,j}^+-q_{Z,j}^+)^2}{2q_{Z,j}^+}+
\frac{(h_{Z,j}^--q_{Z,j}^-)^2}{2q_{Z,j}^-}
+\frac{(h_{Z,j}^++h_{Z,j}^--q_{Z,j}^+-q_{Z,j}^-)^2}{2(2-q_{Z,j}^+-q_{Z,j}^-)}=
\frac{(q_{Z,j}^++q_{Z,j}^-)\Delta_j^2}{\la_j^2(2-q_{Z,j}^+-q_{Z,j}^-)},
\end{eqnarray*}
where
$$
\Delta_j=q_{Z,j}^++q_{Z,j}^--h=h(e^{d_j^+}\1_{\CA_j^+}+
e^{d_j^-}\1_{\CA_j^-}-2)/2.
$$
Since $\max_{1\le j\le p}q_{Z,j}^{\pm}\le 1/2$, we only have  to control the
sum $\sum_{j=1}^p\Delta_j^2$.
\begin{eqnarray*}
E^X_0 \1_{\CA}\Delta_j^2&\leq &\frac{h^2}{2}\Big(e^{a_j^2}\Phi(T_j-2a_j)
+e^{-a_j^2}-4\Phi(T_j-a_j)+2\Big)\ .
\end{eqnarray*}

\noindent
{\bf CASE 1:} $nb^2= O(1)$. By corollary \ref{C1},
$(a_j/(\sqrt{n}b)-1)=o_{P_X}(1)$. Consequently, $\Phi(T_j-2a_j)=
1-o_{P_X}(p^{-2})$.
$$E^X_0
\left[\1_{\CA}\sum_{j=1}^p\Delta_j^2\right]\leq \frac{ph^2}{2}
\sinh^2(nb^2(1+o_{p_X}(1))/2)+o_{P_X}(1)= \frac{ph^2nb^2}{2}+
o_{P_X}(1)=o_{P_X}(1)\ , $$
since $r^2_{n,p}=o(\sqrt{p}/n)$. \\

\noindent
{\bf CASE 2:} $\lim \sup nb^2=\infty$. This implies that
$k^2=o(p)$ and
therefore $ph^2=o(1)$.

\begin{eqnarray*}
E^X_0
\left[\1_{\CA}\sum_{j=1}^p\Delta_j^2\right]=\frac{ph^2}{2}
e^{nb^2(1+o_{P_X}(1))}+o(1)=p^{-(2\beta-1)+ x^2+o_{P_X}(1)} +o(1)\ .
\end{eqnarray*}
Since $x<\varphi(\beta)\leq \sqrt{2\beta-1}$ for $\beta>1/2$, this allows to
conclude.
\endproof

\subsubsection{Study of $E_{\tilde\pi_{Z}}\Delta^2$}\label{SS3}

By Lemma \ref{L0}, the relation \nref{probab*} follows from $
\tilde\pi(\Sigma)\to 1, $ in $P_0$-probability. Thus, we only need to
check that in $P_0$-probability, $ E_{\tilde\pi_{Z}}\Delta^2\to 0$
for $\Delta=\Delta(X,\theta)$ defined by \nref{gD}. By Markov
inequality, the last relation follows from
$$
E_0^X(E_{\tilde\pi_{Z}}\Delta^2)))\to 0,\quad \text{in
$P_X$-probability}.
$$

Let us introduce the events $\CX_{n,p}$. Taking a
positive family $\eta=\eta_{n,p}\to 0$, we set
$$
\CX^j=\{(\|X_j\|^2-n)<\eta n\},\quad
\CX^{ij}=\{\log(p)|(X_j,X_l)|<\eta n\},\quad \CX_{n,p}=\bigcap_{1\le
j<l\le p}\l(\CX^j\cap\CX^{ij}\r).
$$
It follows from Corollary \ref{C1} that, under assumptions {\bf B2} or {\bf B3}
we can take $\eta=\eta_{n,p}\to 0$ such that $P_X(\CX_{n,p})\to 1$.
We have
$$
E_{\tilde\pi_{Z}}\Delta^2=b^4E_{\tilde\pi_{Z}}\l(\sum_{j_1,j_2,j_3,j_4=1}^p
\theta_{j_1}
\theta_{j_2}\theta_{j_3}\theta_{j_4}(X_{j_1},X_{j_2})(X_{j_3},X_{j_4})\r)=
2A_2+6A_3+24A_4,
$$
where
\begin{eqnarray}\label{A2}
A_2&=&b^4\sum_{1\le j_1<j_2\le p}E_{\tilde\pi_{Z}}\l(\epsilon_{j_1}^2
\epsilon_{j_2}^2\r)(X_{j_1},X_{j_2})^2,\\
\label{A3}
A_3&=&b^4\sum_{1\le j_1<j_2<j_3\le
p}E_{\tilde\pi_{Z}}\l(\epsilon_{j_1}^2
\epsilon_{j_2}\epsilon_{j_3}\r)(X_{j_1},X_{j_2})(X_{j_1},X_{j_3}),\\
\label{A4}
A_4&=&b^4\sum_{1\le j_1<j_2<j_3<j_4\le
p}E_{\tilde\pi_{Z}}\l(\epsilon_{j_1}
\epsilon_{j_2}\epsilon_{j_3}\epsilon_{j_4}\r)(X_{j_1},X_{j_2})(X_{j_3},X_{j_4}).
\end{eqnarray}

\subsubsection{Expectation over $\tilde\pi_{Z}$ and over
$E_0^X$}\label{4.3.1}

Let us define the variables $\eta_k$ in $\{1,-1\}$.  The expectations  over
$\tilde\pi_{Z}$ are of the form
\begin{eqnarray}\label{E2}
E_{\tilde\pi_{Z}}\l(\epsilon_{j_1}^2
\epsilon_{j_2}^2\r)&=&\frac{(q_{j_1}^++q_{j_1}^-)(q_{j_2}^++q_{j_2}^-)}{4}=\frac
{1}{4}
\sum_{\eta_1,\eta_2}\prod_{k=1}^2q_{j_k}^{\eta_k},\\
\label{E3} E_{\tilde\pi_{Z}}\l(\e_{j_1}^2
\e_{j_2}\e_{j_3}\r)&=&\frac{(q_{j_1}^++q_{j_1}^-)(q_{j_2}^+-q_{j_2}^-)
(q_{j_3}^+-q_{j_3}^-)}{8}=
\sum_{\eta_1,\eta_2,\eta_3}\frac{\eta_2\eta_3}{8}\prod_{k=1}^3q_{j_k}^{\eta_k}\\
\label{E4} E_{\tilde\pi_{Z}}\l(\e_{j_1}
\e_{j_2}\e_{j_3}\e_{j_4}\r)&=&\frac{1}{16}\prod_{k=1}^4(q_{j_k}^+-q_{j_k}
^-)=\frac{1}{16}
\sum_{\eta_1,\eta_2,\eta_3,\eta_4}\eta_1\eta_2\eta_3\eta_4
\prod_{k=1}^4q_{j_k}^{\eta_k}.
\end{eqnarray}
Let us take the expectation $E^X_0$ over $Y$ of each of these expressions.
We define the vector  $V=b\sum_{k=1}^m\eta_kX_{j_k}$. Here, $E_V^X$ refers
to the expectation of $Y$ over the Gaussian measure $\CN(V,I_n)$. We derive that
\begin{eqnarray*}
&&E^X_0\l(
\prod_{k=1}^mq_{j_k}^{\eta_k}\r)=\frac{h^m}{2^m}E^X_0\l(e^{-\frac
12\sum_{k=1}^mb^2\|X_{j_k}\|^2+b(Y,\sum_{k=1}^m\eta_kX_{j_k})}
\prod_{k=1}^m\1_{\{(Y,\eta_kX_{j_k})<T_{j_k}\|X_{j_k}\|\}}\r)\\
&=&\frac{h^m}{2^m}\exp\l(b^2{\sum_{1\le r<s\le
m}\eta_r\eta_s(X_{j_r},X_{j_s})}\r)E^X_0\l(e^{-\frac
12\|V\|^2+(Y,V)}\prod_{k=1}^m\1_{\{(Y,\eta_kX_{j_k})<T_{j_k}\|X_{j_k}\|\}}\r)\\
&=&\frac{h^m}{2^m}\exp\l({b^2\sum_{1\le r<s\le
m}\eta_r\eta_s(X_{j_r},X_{j_s})}\r)E^X_V\l(\prod_{k=1}^m\1_{\{(Y,\eta_kX_{j_k})
<T_{j_k}\|X_{j_k}\|\}}\r)\\
&=&\frac{h^m}{2^m}\exp\l({b^2\sum_{1\le r<s\le
m}\eta_r\eta_s(X_{j_r},X_{j_s})}\r)P_{j_1,...,j_m}(\eta),
\end{eqnarray*}
where
\begin{eqnarray*}
P_{j_1,...,j_m}(\eta)&=&E^X_V\l(\prod_{k=1}^m\1_{\{(Y,\eta_kX_{j_k})<T_{j_k
}\|X_{j_k}\|\}}\r)
\\
&=&E^X_0\l(\prod_{k=1}^m\1_{\{(Y+V,\eta_kX_{j_k})<T_{j_k}\|X_{j_k}\|\}}\r)\\
&=&
E^X_0\l(\prod_{k=1}^m\1_{\{(Y,\eta_kX_{j_k})<
T_{j_k}\|X_{j_k}\|-(V,\eta_kX_{j_k})\}}\r)\\
&=&E^X_0\l(\prod_{k=1}^m\1_{\{\eta_ky'_{j_k}<T_{j_k}-
(V,\eta_kX_{j_k})/\|X_{j_k}\|\}}\r)\ .
\end{eqnarray*}
Let us define
$$
m_{j_k}(\eta)=\eta_k\sum_{s=1,\,s\not=k}^m\eta_s(X_{j_s},X_{j_k})/\|X_{j_k}\|,
\quad  z_{k}=\eta_k y'_{j_k}.
$$
Then, $P_{j_1,...,j_m}(\eta)$ writes as
$$P_{j_1,...,j_m}(\eta)=
P^X_0\l(z_1<T_{j_1}-a_{j_1}-bm_{j_1}(\eta) \ldots,z_m<
T_{j_m}-a_{j_m}-bm_{j_m}(\eta)\r) .$$
We have
\begin{eqnarray}\label{exp_2}
E^X_0\l( \prod_{k=1}^2q_{j_k}^{\e_k}\r)&=&\frac{h^2}{4}
\exp\l(\eta_1\eta_2 b^2(X_{j_1},X_{j_2})\r)P_{j_1,j_2}(\eta) ,\\
\label{exp_3} E^X_0\l( \prod_{k=1}^3q_{j_k}^{\e_k}\r)&=&\frac
{h^3}{8} \exp\l(b^2\sum_{1\le s<r\le
3}\eta_{s}\eta_{r}(X_{j_s},X_{j_r})\r)P_{j_1,j_2,j_3}(\eta),\\
\label{exp_4} E^X_0\l(\prod_{k=1}^4q_{j_k}^{\e_k}\r)&=& \frac
{h^4}{16}\exp\l(b^2\sum_{1\le s<r\le
4}\eta_{s}\eta_{r}(X_{j_s},X_{j_r})\r)P_{j_1,j_2,j_3,j_4}(\eta).
\end{eqnarray}

\subsubsection{Evaluation of probabilities
$P_{j_1,...,j_m}(\eta)$}\label{EvProb}
By definition of $(z_1,\ldots, z_m)$ we have
$$
E^X_0z_k=0,\quad E^X_0z_k^2=1,\quad E^X_0z_kz_s\eq
r_{ks}(\eta)=\frac{\eta_k\eta_s(X_{j_k},X_{j_s})}{\|X_{j_k}\|\|X_{j_s}\|}
,\quad 1\le k<s\le m.
$$
Denote $\tilde T_{j_k}=T_{j_k}-a_{j_k}$.
Observe that
\begin{eqnarray}\nonumber
P_{j_1,...,j_m}(\eta)&=&1-\sum_{k=1}^m\Phi(-\tilde T_{j_k}-bm_{j_k}(\eta))\\
&+&O\l(\sum_{1\le k<s\le m}P_{r_{ks}(\eta)}\l(- \tilde T
_{j_k}-bm_{j_k}(\eta),- \tilde
T_{j_s}-bm_{j_s}(\eta)\r)\r),\qquad\label{REL}
\end{eqnarray}
where we set, for the Gaussian random vector $(z_1,z_2)$ with
$Ez_k=0, Ez_k^2=1, k=1,2, Ez_1z_2=r$,
$$
P_r(t_1,t_2)=P(z_1<t_1,z_2<t_2)=P(z_1>-t_1,z_2>-t_2).
$$

\noindent The control of $P_{j_1,...,j_m}(\eta)$ then depends on the sequence
$x_{np}$.

\noindent
{\bf CASE 1:} $x=0$. Under the event $\CX_{n,p}$, we have
$\max_{j}a_j=o(\sqrt{\log(p)})$ and  $\tilde T_{j_k}/\sqrt{\log(p)}
\to \infty$. Under the event $\CX_{np}$, we have
$$b|m_{jk}(\eta)|\leq b\sum_{s\neq
k}|(X_{j_s},X_{j_k})|/\|X_{j_k}\|\leq
o(b\sqrt{n}/\log(p))=o(1/\sqrt{\log(p)})\ .$$ It
follows that
$$\max_{j}\Phi(-\tilde
T_{j_k}-bm_{j_k}(\eta))=o(p^{-\alpha}),\quad \forall\ \alpha>0.$$
We conclude
\begin{equation}\label{REL1}
P_{j_1,...,j_m}(\eta)=1-O\l(\sum_{k=1}^m\Phi(-\tilde
T_{j_k}-bm_{j_k}(\eta))\r)=1-o(p^{-\alpha}),\quad \forall\ \alpha>0.
\end{equation}

\noindent
{\bf CASE 2:} $x>0$. We have under the event $\CX_{n,p}$,
$b|m_{jk}(\eta)|=o(b\sqrt{n}/\log(p))=o(1)$ and $\tilde
T_{jk}b= O(\log(p)/\sqrt{n})$. Hence, $\tilde
T_{jk}b|m_{jk}(\eta)|=o(1)$.
Applying Lemma \ref{L4}, we
bound the first
term in \nref{REL}
\begin{eqnarray*}\nonumber
\sum_{k=1}^m\Phi(- \tilde
T_{j_k}-bm_{j_k}(\eta))&=&\sum_{k=1}^m\Phi(- \tilde T_{j_k})-b\sum_{k=1}^m
m_{j_k}(\eta)\Phi( -\tilde
T_{j_k})\\&+&b^2\sum_{k=1}^mO\l(m_{j_k}^2(\eta) \tilde T_{j_k}\Phi(
-\tilde T_{j_k}) \r)\ .
\end{eqnarray*}
Let us define
\begin{eqnarray}\label{c_{ks}1}
R_m\eq\sum_{k=1}^m\Phi(- \tilde T_{j_k})=o((ph)^{-1}),
\end{eqnarray}
by \nref{Phi0} and \nref{Phi1} since $x<\varphi(\beta)\leq
\sqrt{2}(1-\sqrt{1-\beta})$. Applying again \nref{Phi0} and \nref{Phi1}, we get
\begin{eqnarray*}
b(X_{j_s},X_{j_k})\Phi(- \tilde T_{j_k})/\|X_{j_k}\|&=&O\left[
T_{j_k}(X_{j_s},X_{j_k})\Phi(- \tilde
T_{j_k})n^{-1}\right]\\ &=& O\left[
T_{j_k}\Phi( - \tilde T_{j_k})|r_{ks}|\right]=o(|r_{ks}|/(ph))\\
b^2(X_{j_s},X_{j_k})^2\tilde T_{j_k}\Phi( -\tilde
T_{j_k})n^{-1}&=& O(T_{j_k}^3)\Phi( -\tilde
T_{j_k})r_{ks}^2=o(r_{ks}^2/(ph)).
\end{eqnarray*}
It follows that
\begin{eqnarray*}
 \sum_{k=1}^m\Phi(- \tilde
T_{j_k}-bm_{j_k}(\eta))&=&R_m -\sum_{1\le k<s\le
m}o\l(|r_{ks}|/(ph)\r)+\sum_{1\le k<s\le
m}o\l(r_{ks}^2/(ph)\r)\\
&=&\sum_{1\le k<s\le
m}o\l((1+|r_{ks}|+r_{ks}^2)/(ph)\r).
\end{eqnarray*}
Let us turn to the second term in \nref{REL}. If $\Tilde{T}_{j_k}\geq \log(p)$,
then
\begin{eqnarray*}
P_{r_{ks}(\eta)}\l(- \tilde T _{j_k}-bm_{j_k}(\eta),- \tilde
T_{j_s}-bm_{j_s}(\eta)\r)\leq \Phi(- \tilde T _{j_k}-bm_{j_k}(\eta)) =
o(\l((ph)^{-2}\r)
\end{eqnarray*}
If $\tilde{T}_{j_k}\leq \log(p)$, we have
$\tilde{T}_{j_k}r_{ks}=o(1)$ under the event $\CX_{np}$.
By Lemma \ref{L3}  and previous evaluations, we get
\begin{eqnarray*}
&&P_{r_{ks}(\eta)}\l(- \tilde T _{j_k}-bm_{j_k}(\eta),- \tilde
T_{j_s}-bm_{j_s}(\eta)\r)\\&&=\Phi(- \tilde
T_{j_k}-bm_{j_k}(\eta))\,\Phi(- \tilde
T_{j_s}-bm_{j_s}(\eta))O\l(1+r_{ks}^2+|r_{ks}|\r)
=o\l((ph)^{-2}\r).
\end{eqnarray*}
Finally, we obtain
\begin{eqnarray}
P_{j_1,...,j_m}(\eta)&=&1-R_m+o\l(\sum_{1\le k<s\le
m}|r_{ks}|/ph\r)+o((ph)^{-2}) \label{prob*} \\
&=&1+o((ph)^{-1}). \label{prob}
\end{eqnarray}

\subsubsection{Evaluation of $A_2$}

We have $b^2\max_{1\le j_1<j_2\le p}|(X_{j_1},X_{j_2})|=o(1)$ under
the event $\CX_{n,p}$. Since $P_{j_1,j_2}(\eta)=O(1)$, we get from
\nref{exp_2}
$$
E^X_0\l(\prod_{k=1}^2q_{j_k}^{\e_k}\r)=O(h^2).
$$
By Assumption ${\bf B1}$, we have
$$
\sup_{j_1\neq j_2}\  E_X(X_{j_1},X_{j_2})^2=O(n)\ .
$$
It then follows from  \nref{A2} and \nref{E2} that $A_2$ is of the order
$$
b^4h^2\sum_{1\le j_1<j_2\le p}(X_{j_1},X_{j_2})^2\asymp
p^2h^2b^4n\asymp nk^2b^4\to 0\ ,
$$
in $P_X$-probability.
\subsubsection{Evaluation of $A_3$}

Let us turn to $A_3$.
 Consider $\eta_k$ as independent random variables taking
values in $\{-1,1\}$ with probabilities $1/2$. By \nref{E3} and
\nref{exp_3}, we can write
\begin{eqnarray*}
E_{\tilde\pi_{Z}}\l(\theta_{j_1}^2
\theta_{j_2}\theta_{j_3}\r)=\frac{h^3}{8}E_{\eta}\l(\eta_2\eta_3\exp\l(b^2\sum_{
1\le
s<r\le
3}\eta_{s}\eta_{r}(X_{j_s},X_{j_r})\r)P_{j_1,j_2,j_3}(\eta)\r).
\end{eqnarray*}
 Under
the event $\CX_{n,p}$ it follows from \nref{REL1}, \nref{prob}, and the
definition of $\CX_{n,p}$ that
$$
 P_{j_1,j_2,j_3}(\eta)=1+o((ph)^{-1})\text{ and }
b^2\sum_{1\le s<r\le 3}\eta_{s}\eta_{r}(X_{j_s},X_{j_r})=o(1).
$$
It follows that
\begin{eqnarray*}
E_{\tilde\pi_{Z}}\l(\theta_{j_1}^2
\theta_{j_2}\theta_{j_3}\r)&=&\frac{h^3}{8}E_{\eta}
\l(\eta_2\eta_3\exp\l(b^2\sum_ { 1\le
s<r\le 3}\eta_{s}\eta_{r}(X_{j_s},X_{j_r})\r)\r)+o\l(h^2p^{-1}\r).
\end{eqnarray*}
By Taylor expansion of the exponential function, the expectation over $\eta$ is
of the form, for $c_{sr}=b^2(X_{j_s},X_{j_r})$,
\begin{eqnarray*}
E_{\eta}\l(\eta_2\eta_3\l(1+\eta_1\eta_2 c_{12}+\eta_1\eta_3
c_{13}+\eta_2\eta_3 c_{23}
+O\l(c^2_{12}+c^2_{13}+c^2_{23}\r)\r)\r)\\
=b^2(X_{j_2},X_{j_3})+O\l(b^4\sum_{1\le s<r\le
3}(X_{j_s},X_{j_r})^2\r).
\end{eqnarray*}
Under the event $\CX_{n,p}$, we derive from \nref{A3} that
$$
A_3\le h^3\l(b^6O(H_1)+b^8O(H_2)\r)+b^4o(H_3h^2p^{-1}),
$$
where
\begin{eqnarray*}
H_1&=&\sum_{1\le j_1<j_2<j_3\le
p}|(X_{j_1},X_{j_2})||(X_{j_1},X_{j_3})||(X_{j_2},X_{j_3})|,\\
H_2&=&\sum_{1\le j_1<j_2<j_3\le
p}\sum_{1\leq
s<r\leq 3}|(X_{j_1},X_{j_2})||(X_{j_1},X_{j_3})|(X_{j_s},X_{j_r})^2,\\
H_3&=&\sum_{1\le j_1<j_2<j_3\le
p}|(X_{j_1},X_{j_2})||(X_{j_1},X_{j_3})|.
\end{eqnarray*}

Since
\begin{eqnarray*}
|(X_{j_1},X_{j_2})||(X_{j_1},X_{j_3})||(X_{j_r},X_{j_s})|&\le&
|(X_{j_1},X_{j_2})|^3+|(X_{j_1},X_{j_3})|^3+|(X_{j_r},X_{j_s})|^3
,\\
|(X_{j_1},X_{j_2})||(X_{j_1},X_{j_3})|(X_{j_2},X_{j_3})^2&\le&
(X_{j_1},X_{j_2})^4+(X_{j_1},X_{j_3})^4+(X_{j_2},X_{j_3})^4,\\
|(X_{j_1},X_{j_2})||(X_{j_1},X_{j_3})|&\le&
(X_{j_1},X_{j_2})^2+(X_{j_1},X_{j_3})^2,
\end{eqnarray*}
we derive from \nref{moments1}
\begin{eqnarray*}
E_XH_1=O(p^3n^{3/2}),\quad E_XH_2=O(p^3n^{2}),\quad E_XH_1=O(p^3n).
\end{eqnarray*}
Applying Markov's inequality yields
$$
H_1=O_{P_X}(p^3n^{3/2}),\quad H_2=O_{P_X}(p^3n^{2}),\quad H_1=O_{P_X}(p^3n).
$$
Combining these bounds, we obtain
$$
A_3= O_{P_X}((b^6h^3p^3n^{3/2})+O_{P_X}(b^8h^3p^3n^{2})+b^4o_{P_X}(b^4h^2p^2n).
$$
Since $b^4k^2n=o(1),\ hp\asymp k,\ b=o(1)$,
we get $A_3=o_{P_X}(1)$.

\subsubsection{Evaluation of $A_4$} Let us evaluate the item $A_4$.
Similarly to $A_3$, we can write
\begin{eqnarray}\label{equation_Eeta}
E_{\tilde\pi_{Z}}\l(\theta_{j_1}
\theta_{j_2}\theta_{j_3}\theta_{j_4}\r)=h^4E_{\eta}\l(\eta_1\eta_2\eta_3\eta_4
\exp\l(b^2\sum_{1\le s<r\le
4}\eta_{s}\eta_{r}(X_{j_s},X_{j_r})\r)P_{j_1,j_2,j_3,j_4}(\eta)\r).
\end{eqnarray}
Under the event $\CX_{n,p}$ we have
$$
b^2\sum_{1\le s<r\le 4}\eta_{s}\eta_{r}(X_{j_s},X_{j_r})=o(1).
$$

\noindent
{\bf CASE 1:} $x>0$. By \nref{prob*}, we have
$$
P_{j_1,j_2,j_3,j_4}(\eta)=1-R_4+o\l(\sum_{1\le s<r\le
4}|r_{sr}|/hp\r)+o((ph)^{-2}).
$$
Applying a Taylor expansion of the
exponential term in \nref{equation_Eeta} yields
\begin{eqnarray*}
\lefteqn{E_{\eta}\l(\eta_1\eta_2\eta_3\eta_4
\exp\l(b^2\sum_{1\le s<r\le
4}\eta_{s}\eta_{r}(X_{j_s},X_{j_r})\r)P_{j_1,j_2,j_3,j_4}(\eta)\r) }&&\\
&=& E_{\eta}\l(\eta_1\eta_2\eta_3\eta_4\l(1+b^2\sum_{1\le s<r\le
4}\eta_{s}\eta_{r}(X_{j_s},X_{j_r})\r)(1-R_4)\r)+O(|\delta_1|)+O(|\delta_2|)\\
&= & O(|\delta_1|)+O(|\delta_2|),
\end{eqnarray*}
where
\begin{eqnarray*}
\delta_1&=&O\l(b^4\sum_{1\le
s<r\le 4}(X_{j_s},X_{j_r})^2\r),\\
\delta_2&=&o\l(\sum_{1\le k<s\le 4}|r_{ks}|/ph\r)+o((ph)^{-2})\ .
\end{eqnarray*}
{\bf CASE 2:} $x=0$. By \nref{REL1},
$P_{j_1,j_2,j_3,j_4}(\eta)=1-o(p^{-2})$. Arguing as in Case 1, we get
\begin{eqnarray*}
\lefteqn{E_{\eta}\l(\eta_1\eta_2\eta_3\eta_4
\exp\l(b^2\sum_{1\le s<r\le
4}\eta_{s}\eta_{r}(X_{j_s},X_{j_r})\r)P_{j_1,j_2,j_3,j_4}(\eta)\r) }&&\\
&=& O\l(b^4\sum_{1\le
s<r\le 4}(X_{j_s},X_{j_r})^2\r) + o(p^{-2}).
\end{eqnarray*}
All in all, we obtain that under the event $\CX_{n,p}$,
$$
A_4\le h^4b^8O(H_1)+\begin{cases}o(H_2b^4h^4/p^2),&x=0,\\
o(H_3b^4h^3/np+H_2b^4h^2/p^2),&x>0,\end{cases}
$$
where
\begin{eqnarray*}
H_1&=&\sum_{1\le j_1<j_2<j_3<j_4\le
p}|(X_{j_1},X_{j_2})||(X_{j_3},X_{j_4})|\sum_{1\le s<r\le 4}(X_{j_s},X_{j_r})^2,
\\
H_2&=&\sum_{1\le j_1<j_2<j_3<j_4\le
p}|(X_{j_1},X_{j_2})||(X_{j_3},X_{j_4})|,\\
H_3&=&\sum_{1\le j_1<j_2<j_3<j_4\le
p}|(X_{j_1},X_{j_2})||(X_{j_3},X_{j_4})|\sum_{1\le s<r\le
4}|(X_{j_s},X_{j_r})|.
\end{eqnarray*}
We combine the classical upper bounds,
\begin{eqnarray*}
|(X_{j_1},X_{j_2})||(X_{j_3},X_{j_4})|(X_{j_s},X_{j_r})^2&\le&
(X_{j_1},X_{j_2})^4+(X_{j_3},X_{j_4})^4+(X_{j_s},X_{j_r})^4,\\
|(X_{j_1},X_{j_2})||(X_{j_3},X_{j_4})|&\le&
(X_{j_1},X_{j_2})^2+(X_{j_3},X_{j_4})^2,\\
|(X_{j_1},X_{j_2})||(X_{j_3},X_{j_4})||(X_{j_s},X_{j_r})|&\le&
|(X_{j_1},X_{j_2})|^3+|(X_{j_3},X_{j_4})|^3+|(X_{j_s},X_{j_r})|^3\ .
\end{eqnarray*}
with  \nref{moments1} and obtain
\begin{eqnarray*}
E_X(H_1) =O(p^4n^2),\quad E_X(H_2) =O(p^4n),\quad E_X(H_3)
=O(p^4n^{3/2}).
\end{eqnarray*}
Applying Markov's inequality yields
$$
H_1=O_{P_X}(p^4n^2),\ H_2=O_{P_X}(p^4n),\ H_3=O_{P_X}(p^4n^{3/2}).
$$
Since $b^4k^2n=o(1),\ hp\asymp k$, we get
$$
h^4b^8H_1=O(h^4p^4b^8n^2)=o_{P_X}(1),\quad H_2b^4h^2/p^2=O_{P_X}(b^4h^2p^2n)=
o_{P_X}(1).
$$
If $x>0$, we also have to upper bound the term $H_3$. Since
$r_{n,p}^2=o(1/\sqrt{n})$ (cf. (\ref{L1})) and since $x>0$, we
derive that $k=o(\sqrt{n})$. Then, we get
$$
H_3b^4h^3/np=O_{P_X}(b^6p^3h^3n^{3/2}/nb^2)=o_{P_X}(1),\quad
$$
Therefore we obtain $A_4=o_{P_X}(1)$. The
proposition follows. \endproof

\subsection{Proof of Proposition \ref{P2}}\label{ProofP2}

We will prove that there exists a family of events $\CZ_{n,p}$
such that $P_0(\CZ_{n,p})\to 1$ and
$$
\log(\Lambda(Z))=\sum_{j=1}^p\log\l(1+(h/2)\l(e^{d_j^+}+e^{d_j^-}-2\r)\r)\to
0,\quad Z\in \CZ_{n,p}.
$$
We take $\CZ_{n,p}=\{(X,Y): |y'_j|\le T_j,\ 1\le j\le p,\
X\in\CX_{n,p}\}$ where $\CX_{n,p}$ was defined in Section \ref{SS3}.
It follows from Lemma \ref{L0} and Section \ref{SS3} that
$P_0(\CZ_{n,p})\to 1$.

Under the events $\CZ_{n,p}$ we can replace the quantities
$(h/2)e^{d_j^\pm}/2$ by $q_j^{\pm}=(h/2)e^{d_j^\pm}\1_{\pm
y'_j<T_j}$, cf. Section \ref{SS2}. Let us consider
$$
\tilde L=\sum_{j=1}^p\log(1+\Delta_j),\quad
\Delta_j=(q_j^++q_j^--h).
$$
Under the event
$\CA=\CA_{n,p}=\bigcap_{j=1}^p(\CA_j^{+}+\CA_j^{-})$ defined in Section
\ref{SS2}, we have uniformly in $1\le j\le p$,
\begin{eqnarray*}
q_{j}^{+}+q_{j}^{-}&=&\frac{h}{2}e^{-a_j^2/2}\cosh(a_jy'_j)\le
he^{-a_j^2/2}\cosh(a_jT_j)\\
&=&(1+e^{-2a_jT_j})/2\sim 1/2,
\end{eqnarray*}
as $h\to 0$. Consequently, we have
$$
\tilde L=\sum_{j=1}^p=A_1+O(A_2),\quad
A_1=\sum_{j=1}^p\Delta_j,\quad A_2=\sum_{j=1}^p\Delta_j^2.
$$
Thus, we need to show that $A_1\to 0$ and that $A_2\to 0$ in
$P_0$-probability. It was stated in the proof of Lemma
\ref{Lem2} that $E_0^XA_2=o_{P_X}(1)$.  Markov's inequality then allows to
derive that $A_2=o_{P_X}(1)$.  In order to prove the first relation, we shall
show
that $E_0^XA_1\to 0$ and that $\Var_0^XA_1\to 0$ in $P_0$-probability.
Observe that
\begin{eqnarray*}
E_0^XA_1=h\sum_{j=1}^p\l(\Phi(T_j-a_j)-1\r) =
-h\sum_{j=1}^p\Phi(-T_j+a_j).
\end{eqnarray*}
By \nref{Phi0} and \nref{Phi1} we have
$$
h\sum_{j=1}^p\Phi(-T_j+a_j)\asymp \sum_{j=1}^p\Phi(-T_j)=o(1).
$$
We have $\Var_0^XA_1\le B+A_2$ with  $B=\sum_{1\le j<l\le
p}\hat\Delta_j\hat\Delta_l$ and  $\hat\Delta_j=\Delta_j-E_0^X\Delta_j$. We need
to check that, in $P_X$ probability,
$$
E^X_0(B)=\sum_{1\le j<l\le p}E^X_0(\hat\Delta_j\hat\Delta_l)\to 0.
$$
Note that
$$
E^X_0(\hat\Delta_j\hat\Delta_l)=B_{jl}-C_{jl},
$$
where
$$
B_{jl}=E^X_0\l((q_j^++q_j^-)(q_l^++q_l^-)\r),\quad
C_{jl}=h^2\Phi(T_j-a_j)\Phi(T_l-a_l)\ .
$$
We consider independent random variables $\eta_1,\eta_2$ taking
values $-1$ and 1 with probabilities $1/2$. We write (compare with
\nref{exp_2})
\begin{eqnarray*}
B_{jl}=h^2E_{\eta}\l[\exp\l({\eta_1\eta_2b^2(X_j,X_l)}\r)P_{j,l}
(\eta)\r]\ ,\quad
C_{jl}=h^2P_{j,l}^0\ .
\end{eqnarray*}
Here we set
$$
P_{j,l}^0=\Phi(\tilde T_j)\Phi(\tilde T_l)=1-\Phi(-\tilde
T_j)-\Phi(-\tilde T_l)+\Phi(-\tilde T_j) \Phi(-\tilde T_l),\ \tilde
T_l=T_l-a_l .
$$
We obtain the new decomposition
\begin{equation}\label{item}
B_{jl}-C_{jl}=h^2(U_{jl}+V_{jl})\ ,
\end{equation}
where
$$
U_{jl}=E_{\eta}\l[\l(\exp\l(({\eta_1\eta_2b^2(X_j,X_l)}\r)-1\r)P_{j,l}(\eta)\r],
\quad
V_{jl}=E_{\eta}\l(P_{j,l}(\eta)-P_{j,l}^0\r)\ .
$$
Let us recall some notations introduced in Section \ref{EvProb}.
$r_{jl}(\eta)=\eta_1\eta_2r_{jl}$,
$$
r_{jl}=\frac{(X_j,X_l)}{\|X_j\|\|X_l\|},\quad
m_{jl}(\eta)=\frac{\eta_1\eta_2(X_j,X_l)}{\|X_j\|},\quad
m_{lj}(\eta)=\frac{\eta_1\eta_2(X_j,X_l)}{\|X_l\|}\ .
$$
Moreover,  $z_j$
and $z_l$ stand for  standard Gaussian variables  with
$\Cov(z_j,z_l)=r_{jl}(\eta)$. Then, $P_{j,l}(\eta)$ is written as
\begin{eqnarray*}
P_{j,l}(\eta)&=& P_0^X(z_j<\tilde T_j-bm_{jl}(\eta),\,z_l<\tilde T_l-
bm_{lj}(\eta))\\
&=&1-\Phi(-\tilde T_j+bm_{jl}(\eta))-\Phi(-\tilde T_l+bm_{lj}(\eta))\\
&&+ P_0^X(z_j<-\tilde T_j+bm_{jl}(\eta),\,z_l<-\tilde T_l+
bm_{lj}(\eta))\ .
\end{eqnarray*}

\noindent {\bf CASE 1:} $x=0$. The evaluations of the terms $V_{jl}$ in
\nref{item} are  similar to the
ones in Section \ref{EvProb}. We get
$$
P_{j,l}^0=1-o(p^{-2}),\quad P_{j,l}(\eta)=1-o(p^{-2}),\quad
|P_{j,l}(\eta)-P_{j,l}^0|=o(p^{-2}).
$$
We derive that $h^2\sum_{1\le j<l\leq p}V_{j,l}=o(h^2)$.\\

\noindent {\bf CASE 2:} $x>0$. We have (compare with
\nref{c_{ks}1} and \nref{prob})
\begin{eqnarray*}
\Phi(-\tilde T_j) \Phi(-\tilde T_l)&=&o((ph)^{-2}),\\
P_0^X(z_j<-\tilde T_j+bm_{jl}(\eta),\,z_l<-\tilde T_l+
bm_{lj}(\eta))&=&o((ph)^{-2}),\\
 \Phi(-\tilde T_j+bm_{jl}(\eta))&=&\Phi(-\tilde T_j)+\eta_1\eta_2
br_{jl}+o(r_{jl}^2/(ph))\ ,
\end{eqnarray*}
Taking the
expectation over $\eta$, we get
$$
E_\eta \l(P_{j,l}(\eta)\r)-P_{j,l}^0=o(r_{jl}^2/(hp))+o((ph)^{-2}).
$$
in $P_X$-probability. Therefore
$$
h^2\sum_{1\le j<l\le p}V_{jl}=O(Hhp^{-1})+o(1),\quad H=\sum_{1\le
j<l\le p}r_{jl}^2 .
$$
 Under $\CX_{n,p}$ we have $r_{jl}^2\sim n^{-2}(X_j,X_l)^2$.
Since $E_X[(X_j,X_l)^2]=O(n)$  for $j\not= l$ (Assumption ${\bf
B1}$), we get
$$
H\sim n^{-2}\sum_{1\le j<l\le p}(X_j,X_l)^2=O_{P_X}(n^{-1}p^2)\ .
$$
This leads to $h^2\sum_{1\le j<l\le p}V_{jl}=
O_{P_X}(ph/n)+o(1)$. Since
$r_{n,p}^2=o(1/\sqrt{n})$ (Eq. \ref{L1}) and since $x>0$, we derive that
$k=o(\sqrt{n})$. Consequently, we have $ph/n=O(k/n)=o(1)$.\\

Let us turn to the terms $U_{jl}$. They are handled as in
 Section
\ref{SS3}. We have
\begin{eqnarray*}
U_{jl}&=&
E_{\eta}\l(\l(\eta_1\eta_2b^2(X_j,X_l)+O\r(b^4(X_j,X_l)^2\r)
\l(1+o((ph)^{-1})\r)\\
&=&O\l(b^4(X_j,X_l)^2\r)+O\l(b^2|(X_j,X_l)|/(ph)\r).
\end{eqnarray*}
Then, we get
$$
h^2\sum_{1\le j<l\le p}U_{jl}= O\l(h^2b^4H_1\r)+O\l(hb^2H_2/p\r),
$$
where
$$
H_1=\sum_{1\le j<l\le p}(X_j,X_l)^2,\quad H_2= \sum_{1\le j<l\le
p}|(X_j,X_l)|\le pH_1^{1/2}.
$$
Arguing as for $H$,  we get
$$
H_1=O_{P_X}(p^2n),\quad H_2=O_{P_X}(p^2n^{1/2}).
$$
It follows that
$$
\sum_{1\le j<l\le
p}B_{jl}=O_{P_X}(p^2h^2b^4n)+O_{P_X}(phb^2n^{1/2})=o_{P_X}(1),
$$
since $p^2h^2b^4n\asymp k^2b^4n\to 0$ by \nref{L1}. The proposition follows.
\endproof

\subsection{Proof of Theorem \ref{prte_lower_bound_unknown}}

An in the proof of Theorem \ref{TL}, we consider $x=\lim \sup x_{n,p}$ and
we take  $c\in(0,1)$ such that $x_c=x/c<\varphi(\beta)$. We also
define $b=x_c\sqrt{\log(p)/n}$. We first consider the case where $k\log(p)/n
\rightarrow 0$.

We use a different prior $\pi$ than for Theorem \ref{TL}. Let us note
$\mathcal{M}(k,p)$ the collection of subsets
of $\{1,\ldots, p\}$ of size $k$. We consider a random
vector $\theta=(\theta_j)$ with coordinates $\theta_j=b\epsilon_j$ where
$\epsilon_j\in(0,1)$. The set of non-zero coefficient of $\epsilon$ is drawn
uniformly in $\mathcal{M}(k,p)$. This introduces a prior probability $\pi$ on
$\theta$.

Consider the mixture
$$P_{\pi}(dZ)= E_{\pi}P_{\theta,\sqrt{1- bk^2 }}(dZ)=
\int_{\mathbb{R}^p}P_{\theta,\sqrt{1- bk^2 }}(dZ)\pi(d\theta)$$
and the likelihood ratio
$$L_{\pi}(Z)= \frac{dP_{\pi}}{dP_{0,1}}(Z)\ .$$
As in the proof of Theorem  \ref{TL}, we shall prove that $L_{\pi}(Z)$
converges to $1$ in $P_0$ probability. This will enforce that
$\gamma^{un}_{n,p,k}[x_c\sqrt{k\log(p)}/\sqrt{1-kb^2}]\rightarrow 1$. Since
$kb^2$ converges to $0$, this will complete the proof.\\

\noindent The likelihood ratio has the form $L_{\pi}(Z)=
\sum_{m\in\mathcal{M}(k,p)}|\mathcal{M}(k,p)|^{-1}L_{m}(Z)$ and
\begin{eqnarray}
 L_{m}(Z)&=&\nonumber (1-kb^2)^{-n/2}
\exp\left(-\frac{kb^2\|Y\|^2}{2(1-kb^2)}+ \frac{b(Y,\sum_{i\in m}X_i)}{1-kb^2}
\right)\\ &\times & \exp\left[-\sum_{i,j\in m}
\frac{b^2}{2(1-kb^2)}(X_i,X_j)\right]\ .\label{ratiolikelihood}
\end{eqnarray}

\begin{defi}
Consider $\delta\in(0,1)$, a positive integer $s$ and  a $n\times p$ matrix $A$.
 We say that $A$ satisfies a $\delta$-restricted isometry property of order $s$
if for all $\theta\in\mathbb{R}^p_s$,
$$(1-\delta)\|\theta\|\leq \|A\theta\|\leq (1+\delta)\|\theta\|\ .$$
\end{defi}

Let us define the events $\Omega_1$ and $\Omega_2$ by
$$\Omega_1:\ "X/\sqrt{n}\text{ satisfies a }\delta^{(1)}_{n,p}\text{ restricted
isometry of
order }2k"$$
$$\Omega_2:\ "\text{For any }1\leq i\leq p,\ (Y,X_i/\|X_i\|)\leq
 \sqrt{2\log(p)}(1+\delta^{(2)}_{n,p})"\ ,$$
where $\delta^{(1)}_{n,p}= 16\sqrt{k\log(p)/n}$ and $\delta^{(2)}_{n,p}=
\log^{-1/2}(p)$. Applying a
deviation inequality due to Davidson and Szarek (Theorem 2.13 in
\cite{davidson2001}), we derive that $P_X(\Omega_1^c)=o(1)$. By the Gaussian
concentration inequality, we have $P_0(\Omega_2^c)=o(1)$.
Then, we take $\Omega=\Omega_1\cap\Omega_2$.

\begin{lemma}\label{lemma1_lowerbound_unknown}
We have  $E_{0}\left[L^2_{\pi}(Z)\1_{\Omega}\right]\leq  1+ o(1)$.
\end{lemma}

\begin{lemma}\label{lemma2_lowerbound_unknown}
We have $E_{0}\left[L_{\pi}(Z)\1_{\Omega^c}\right]= o(1)$.
\end{lemma}
Since $E_{0}\left[L_{\pi}(Z)\right]=1$, we get the desired result by
combining these two
lemmas.\\

Let us turn to the case $k\log(p)/n\rightarrow \infty$. We consider $b>0$
defined by $$\frac{kb^2}{1-kb^2}=(2\beta-1)\frac{k\log(p)}{n}\ .$$
\begin{lemma}\label{lemma3_lowerbound_unknown}
 We have
$$E_{0}\left[L_{\pi}^2(Z)\right]= 1+o(1)\ .$$
\end{lemma}
This lemma implies that for $r= \sqrt{(2\beta-1)k\log(p)/n}\rightarrow \infty$,
we have  $\gamma^{un}_{n,p,k}(r)\rightarrow 1$.
\endproof

In the proof of the following lemmas, $o(1)$ stands for a positive
quantity which depends only on $(k,p,n)$ and tends to $0$ as $(n,p)$
tend to infinity.

\subsubsection{Proof of Lemma \ref{lemma1_lowerbound_unknown}}

In order to
upper bound $E_{0}\left[L^2_{\pi}(Z)\1_{\Omega}\right]$, we first upper bound
$E_{0}\left[L_{m_1}(Z)L_{m_2}(Z)\1_{\Omega}\right]$ for any
$m_1,m_2\in\mathcal{M}(k,p)$. We define $W_1$, $W_2$, $W_3$ by , $W_1=\sum_{i\in
m_1\setminus m_2}X_i$, $W_2=\sum_{i\in m_2\setminus m_1}X_i$, and
$W_3=\sum_{i\in
m_1\cap m_2}X_i$. We note $S=|m_1\cap m_2|$.
\begin{eqnarray*}
L_{m_1}(Z)L_{m_2}(Z)&=&  (1-kb^2)^{-n}
\exp\left(-\frac{kb^2\|Y\|^2}{1-kb^2}+ \frac{b(Y,2W_3+W_1+W_2)}{1-kb^2}
\right)\\ &\times &\exp\left[-
\frac{b^2}{2(1-kb^2)}(\|W_1+W_3\|^2+
\|W_2+W_3\|^2)\right]\ .
\end{eqnarray*}
Let us take the expectation of $L_{m_1}(Z)L_{m_2}(Z)$ with respect to
$(W_1,W_2)$.
$$E_0^{Y,W_3}[L_{m_1}(Z)L_{m_2}(Z)]=
(1-Sb^2)^{-n}\exp\left[-\frac{\|Y\|^2Sb^2}{1-Sb^2}+
\frac{2b(Y,W_3)}{1-Sb^2}- \frac{b^2\|W_3\|^2}{1-Sb^2}\right]\ .$$
When $S=0$, we have $E_0^{Y,W_3}[L_{m_1}(Z)L_{m_2}(Z)]=1$. Let us now consider
the case $S>0$.
On the event
$\Omega$, we have
$$(Y,\frac{W_3}{\|W_3\|})\leq
\sqrt{2\log(p)}(1+\delta^{(2)}_{n,p})\frac{\sum_{i\in m_1\cap
m_2}\|X_i\|}{\|\sum_{i\in m_1\cap m_2}X_i\|} \leq \sqrt{2S\log(p)}
(1+o(1))\ ,$$
since $X/\sqrt{n}$ satisfies a $\delta^{(1)}_{n,p}$-restricted isometry of order
$2k$.
 Then, we can upper bound the expectation with respect to $Y$.
\begin{eqnarray*}
E_0^{W_3}[\1_{\Omega}L_{m_1}(Z)L_{m_2}(Z)]&\leq &
(1-S^2b^4)^{-n/2}\exp\left[\frac{b^2\|W_3\|^2}{1+Sb^2}\right]\\
&\times &\Phi\left[\sqrt{ 2S\log(p)}(1+o(1))-
2b\|W_3\|(1-o(1))\right]\ .
\end{eqnarray*}
Moreover on $\Omega$, we have $\sqrt{1-\delta^{(1)}_{n,p}}\leq
\|W_3\|/\sqrt{nS}\leq
\sqrt{1+\delta^{(1)}_{n,p}}$. Since $k\log(p)/n$ goes to $0$, we get
\begin{eqnarray}
E_0[\1_{\Omega}L_{m_1}(Z)L_{m_2}(Z)]&\leq
&\exp\left[x_c^2S\log(p)(1+o(1))\right]\Phi\left[\sqrt{
S\log(p)}\left(\sqrt{2}-2x_c+o(1))\right)\right]\ .\nonumber
\end{eqnarray}
For any $x<0$, we have $\Phi(x)\leq
e^{-x^2/2}$. Hence, we get  $\Phi(x)\leq
e^{-x_-^2/2}$ for any $x\in\mathbb{R}$. It follows that
\begin{eqnarray}
E_0[\1_{\Omega}L_{m_1}(Z)L_{m_2}(Z)] &\leq & \exp\left[S\log(p)\left\{x_c^2 -
(1-\sqrt{2}x_c)_-^2
+o(1)\right\}\right]\ . \label{majoration_variance_beta}
\end{eqnarray}
Hence, we get
\begin{eqnarray*}
 E_0[\1_{\Omega}L^2_{\pi}(Z)]\leq E_S\left[p^{S\{x_c^2 -
(1-\sqrt{2}x_c)_-^2
+o(1)\}}\right]
\end{eqnarray*}
where $S$ follows a hypergeometric distribution with parameters $p$,
$k$ and $k/p$. We know from Aldous (p.173) \cite{aldous85} that $S$
has the same distribution as the random variable
$E(U|\mathcal{B}_p)$ where $U$ is binomial random variable of
parameters $k$, $k/p$ and $\mathcal{B}_p$ some suitable
$\sigma$-algebra. By a convexity argument, we then obtain
\begin{eqnarray*}
 E_0[\1_{\Omega}L^2_{\pi}(Z)]&\leq& \left[1+\frac{k}{p}\left(p^{x_c^2 -
(1-\sqrt{2}x_c)_-^2
+o(1)}-1\right)\right]^k \\
&\leq & \exp\left[\frac{k^2}{p}p^{x_c^2 -
(1-\sqrt{2}x_c)_-^2
+o(1)}\right]\\
&\leq & \exp\left[p^{1-2\beta+ x_c^2 -
(1-\sqrt{2}x_c)_-^2
+o(1)}\right]
\end{eqnarray*}
Since $x_c<\varphi(\beta)$, one can check that $1-2\beta+ x_c^2 -
(1-\sqrt{2}x_c)_-^2$ is negative and we conclude that $
E_0[\1_{\Omega}L^2_{\pi}(Z)]\leq 1+o(1)$.
\endproof

\subsubsection{Proof of Lemma \ref{lemma2_lowerbound_unknown}}

By symmetry, it is sufficient to prove that $E_0(L_m(Z)\1_{\Omega^c})=o(1)$.
Let us decompose $E_0(L_m(Z)\1_{\Omega^c})= E_0(L_m(Z)\1_{\Omega_2^c})+
E_0(L_m(Z)\1_{\Omega_1^c\cup \Omega_2})$.
Since $E_0^X(L_m(Z))= 1$, $P_X$ almost surely, we have
$E_0(L_m(Z)\1_{\Omega_2^c})= P_X(\Omega_2^c)= o(1)$. Let us turn to
$E_0(L_m(Z)\1_{\Omega_1^c\cup \Omega_2})$. For any $1\leq i\leq p$, we define
the event $\Omega^{(i)}$ by $(Y,X_i/\|X_i\|)\geq
\sqrt{2\log(p)}(1+\delta^{(2)}_{n,p})$.
$$E_0(L_m(Z)\1_{\Omega_1^c\cup \Omega_2})\leq
\sum_{i=1}^pE_0\left[L_m(Z)\1_{\Omega_2}\1_{\Omega^{(i)}}\right]$$
The value of these expectations depends on $i$ through the
property "$i\in m$" or "$i\notin m$". Let us
assume for instance that $1\in m$ and $2\notin m$. Then, we get
\begin{eqnarray}\label{majoration_lemm2_variance_unknown}
 E_0(L_m(Z)\1_{\Omega_1^c\cup \Omega_2})\leq
kE_0\left[L_m(Z)\1_{\Omega_2}\1_{\Omega^{(1)}}\right]+
pE_0\left[L_m(Z)\1_{\Omega_2}\1_{\Omega^{(2)}}\right]\ .
\end{eqnarray}
First, we upper bound $E_0[L_m(Z)\1_{\Omega_2}\1_{\Omega^{(2)}}]$. Taking the
expectation of $L_m(Z)$ with respect to $(X_i)_{i\in m}$ leads to
$E_0^{Y,X_2}[L_m(Z)]=1$. Hence, we get
\begin{eqnarray}\label{majoration2_lemm2_variance_unknown}
E_0[L_m(Z)\1_{\Omega_2}\1_{\Omega^{(2)}}]\leq P_0(\Omega^{(2)})\leq
p^{-1}e^{-\sqrt{\log(p)}}=o(p^{-1})\ .
\end{eqnarray}
Let turn to  $E_0[L_m(Z)\1_{\Omega_2}\1_{\Omega^{(1)}}]$.
We first take the expectation of $L_m(Z)$ conditionally to $X_1$ and $Y$:
$$E_0^{Y,X_1}[L_m(Z)]=
(1-b^2)^{-n/2}\exp\left[-\frac{b^2
\|Y\|^2}{2(1-b^2)}-\frac{b^2\|X_1\|^2}{2(1-b^2)} +
\frac{(Y,X_1)b}{1-b^2} \right ]\ .$$
Then, we take the expectation with respect to $Y$
$$E_0^{X_1}[L_m(Z)\1_{\Omega^{(1)}}]\leq
1-
\Phi\left[\sqrt{\frac{2\log(p)}{1-b^2}}(1+\delta_{n,p}^{(2)})-\frac{
\|X_1\|b } {
\sqrt { 1- b^2} } \right ]
.$$
Moreover, on $\Omega_2$ we have $\|X_1\|\leq \sqrt{n}(1+o(1))$
$$E_0^{X_1}[L_m(Z)\1_{\Omega^{(1)}\cup \Omega_2}]\leq
\Phi\left[\sqrt{\log(p)}(x_c-\sqrt{2}+o(1))\right]\leq
C\exp\left[-\log(p)(\sqrt{2}-x_c-o(1))^2/2\right]$$
for $(n,p)$ large enough, since $x_c<\varphi(\beta)<\sqrt{2}$.
\begin{eqnarray}\label{majoration3_lemm2_variance_unknown}
kE_0^{X_1}[L_m(Z)\1_{\Omega^{(1)}\cup \Omega_2}]\leq
p^{-(\sqrt{2}-x_c)^2/2+ 1-\beta+o(1)}=o(1)\ ,
\end{eqnarray}
since $x_c<\sqrt{2}(1-\sqrt{1-\beta})\leq \varphi(\beta)$. Combining
(\ref{majoration_lemm2_variance_unknown}),
(\ref{majoration2_lemm2_variance_unknown}), and
(\ref{majoration3_lemm2_variance_unknown}) completes the proof.
\endproof

\subsubsection{Proof of Lemma \ref{lemma3_lowerbound_unknown}}
Arguing as in the proof of Lemma \ref{lemma1_lowerbound_unknown}, we get
\begin{eqnarray*}
E_0^{W_3}[L_{m_1}(Z)L_{m_2}(Z)]&= &
(1-S^2b^4)^{-n/2}\exp\left[\frac{b^2\|W_3\|^2}{1+Sb^2}\right]\ .
\end{eqnarray*}
Taking the expectation with respect to $W_3$ leads to
\begin{eqnarray*}
 E_0[L_{m_1}(Z)L_{m_2}(Z)] &= & (1-Sb^2)^{-n/2}\leq
\exp\left[\frac{nSb^2}{2(1-kb^2)}\right]
\end{eqnarray*}
As in the proof  of Lemma \ref{lemma1_lowerbound_unknown}, we upper bound the
term $E_0[L_{\pi}^2(Z)]$ by Jensen's inequality.
\begin{eqnarray*}
E_0[L_{\pi}^2(Z)]&\leq&
\left[1+\frac{k}{p}\left\{\exp\left(\frac{nb^2}{2(1-kb^2)}\right)-1\right\}
\right ] ^k\leq
\exp\left[\frac{k^2}{p}\exp\left(\frac{nb^2}{2(1-kb^2)}\right)\right]\\
&\leq &
\exp\left[p^{1-2\beta}\exp\left\{(\beta-1/2)\log(p)\right\}\right]=1+o(1)\ ,
\end{eqnarray*}
since $b$ satisfies $kb^2/(1-kb^2)= (2\beta-1)k\log(p)/n$.
\endproof

\section{Proofs of the upper bounds}

\subsection{Tests based on the statistic $t_0$}
Recall that
$$
t_0=(2n)^{-1/2}\sum_{i=1}^n(Y_i^2-1).
$$
Under $H_0$, the statistics $Y_i=\xi_i\sim \CN(0,1)$ are i.i.d. This
implies $E_0(t_0)=1,\ \Var_0(t_0)=1$. By the Central Limit Theorem,
$t_0\to \xi\sim \CN(0,1)$ as $n\to\infty$ in $P_0$-probability. This
yields Theorem \ref{TU1} (i).\\

Let us consider the type II errors. We need to show that, if
$nr^4\to\infty$, then $\sup_{\theta\in \Theta_p(r)} P_\t(t_0\le u_\a)\to 0$.
We will prove that, uniformly over $\t\in
\Theta_p(r)$,
\begin{equation}\label{UP1}
E_\t t_0\to\infty,\quad \Var_\t t_0=o((E_\t t_0)^2).
\end{equation}
Indeed, if \nref{UP1} is true, we derive that for
$n,p$ large enough,
\begin{eqnarray}\nonumber
P_\t(t_0\le u_\a)&=&P_\t(E_\t t_0-t_0\ge E_\t t_0- u_\a)\le
P_\t(|E_\t t_0-t_0|\ge E_\t t_0- u_\a)\\
&\le&\frac{\Var_\t(t_0)}{(E_\t t_0- u_\a)^2}=o(1)\ ,\label{cheb}
\end{eqnarray}
by Chebychev's inequality.
In order to check \nref{UP1}, we use the identities
$$
E_\t t_0=E_X(E_\t^X t_0),\quad \Var_\t t_0=\Var_X(E_\t^X
t_0)+E_X(\Var_\t^X t_0).
$$
Under $P_\theta^X,\ \theta \in
\Theta_k(r)$, we have $Y\sim \CN_n(v,I_n)$, where
$$
v=v(\theta,X)=\sum_{j=1}^p\theta_jX_j,\quad
\|v\|^2=\sum_{j=1}^p\theta_j^2\|X_j\|^2+2\sum_{1<j<l\le
p}\theta_j\t_l(X_j,X_l).
$$
It follows that
$$
E_\theta^X(t_0)=(2n)^{-1/2}\|v\|^2, \quad
\Var_\theta^X(t_0)=1+2n^{-1}\|v\|^2. 
$$
Since $E_X(\|X_j\|^2)=n,\ E_X((X_j,X_l))=0,\ j\not=l$, we get the
first convergence in \nref{UP1}:
$$
E_\t
t_0=(2n)^{-1/2}E_X(\|v\|^2)=(n/2)^{1/2}\sum_{j=1}^p\t_j^2=(n/2)^{1/2}\|\t\|^2\ge
(n/2)^{1/2}r^2\to\infty.
$$
Let us turn to the variance term
\begin{eqnarray*}
E_X(\Var_\t^X t_0)&=&1+2n^{-1}E_X(\|v\|^2)=1+2\|\t\|^2=o(E_\t t_0),\\
\Var_X(E_\theta^X(t_0))&=&(2n)^{-1}\Var_X(\|v\|^2).
\end{eqnarray*}
By {\bf A2}, the random variables $X_{ij}$ are independent in
$(i,j), \ i=1,...,n,\ j=1,..,p$. Consequently, the random variables
$(X_{j_1},X_{l_1})$  with $\{j_1,l_1\}\not=\{j_2,l_2\}$ are uncorrelated.
Moreover, $\|X_{j_1}\|^2$ and $(X_{j},X_{l})$ are uncorrelated as long as
$(j,l)\neq (j_1,j_1)$.  We have
$$\Var_X\|X_j\|^2=\Var_X(X_{ij}^2)n,\ E_X(X_j,X_l)^2=n,\
j\not=l,$$ where $ \Var_X(X_{ij}^2)\le E_X(X_{ij}^4)<\infty $ by
{\bf B1}. Then, we get
\begin{eqnarray*}
n^{-1}\Var_X\|v\|^2&=&n^{-1}\sum_{j=1}^p\theta_j^4\Var_X\|X_j\|^2+4n^{-1}
\sum_{1\le
j<l\le p}\theta_j^2\theta_l^2E_X(X_j,X_l)^2\\
&\le& \sup_{i}[E_X(X_{i1}^4)]\sum_{j=1}^p\theta_j^4+4\|\theta\|^4 \le
(O(1)+4)\|\theta\|^4\\
&=&o(n\|\theta\|^4)=o\l((E_\t t_0)^2\r),\quad \text{as}\quad
n\|\theta\|^4\ge nr^4\to\infty.
\end{eqnarray*}
Therefore we get the second relation \nref{UP1}.\\

Note that if $nr^4\to\infty$, then in the inequality
\nref{cheb}, we can replace $u_\a$ by a sequence $T_{np}\to\infty$
such that $\lim\sup T_{np}r^{-2}n^{-1/2}<1$, for instance by
$T_{pn}=n^{1/2}r^2/2$. Then, the corresponding test $\psi^0$ satisfies
$\g(\psi^0,\Theta_p(r))\to 0$.
Theorem \ref{TU1} follows.
\endproof

\subsection{Tests based on the statistic $t_1$}
First observe that under $H_0$, the statistic $t_1$ is a degenerate
$U$-statistic of the second order, i.e., for $Z_s=(X^{(s)},Y_s),\
s=1,2,3$ one has $ E_{Z_1}K(Z_1,Z_2)=0 $, which yields $E_0t_1=0$. By Assumption
{\bf A1},
\begin{eqnarray*}
E_0t_1^2&=&
E_0(K^2(Z_1,Z_2))=p^{-1}E_0(Y_1^2Y_2^2)\sum_{j=1}^p\sum_{l=1}^p
E_X\l(X_{1j}X_{2j}X_{1l}X_{2l}\r)\\
&=&p^{-1}\sum_{j=1}^p E_X\l(X_{1j}^2X_{2j}^2\r)=1.
\end{eqnarray*}
Set
$$
G(Z_1,Z_2)=E_{Z_3} \l(K(Z_1,Z_3)K(Z_2,Z_3)\r),\
G_2=E_0(G^2(Z_1,Z_2)),\ G_4=E_0(K^4(Z_1,Z_2)),
$$
where $E_{Z_3}$ denotes the expectation over ${Z_3}$ under $P_0$. In
order to establish the asymptotic normality of $t_1$ we only need to check the
two following conditions, see \cite{I94} Lemma 3.4,
\begin{equation}\label{Unorm}
G_2=o(1),\quad G_4=o(n^2).
\end{equation}
We have by Assumption {\bf A1},
\begin{eqnarray*}
G(Z_1,Z_2)&=&p^{-1}E_{Z_3}\l(Y_1Y_2Y_3^2\sum_{j=1}^p\sum_{l=1}^p
X_{1j}X_{3j}X_{2l}X_{3l}\r)\\
&=&p^{-1}Y_1Y_2\sum_{j=1}^p\sum_{l=1}^p
X_{1j}X_{2l}E_{X}(X_{3j}X_{3l})\\
&=&p^{-1}Y_1Y_2\sum_{j=1}^pX_{1j}X_{2j}= p^{-1/2}K(Z_1,Z_2).
\end{eqnarray*}
Since $E_0(K^2(Z_1,Z_2))=1$, we get the first convergence in
\nref{Unorm}. Next by {\bf A2},
\begin{eqnarray*}
E_0(K^4(Z_1,Z_2))&=&p^{-2}E_0(Y_1^4Y_2^4)\sum_{j=1}^p\sum_{l=1}^p\sum_{r=1}^p
\sum_{s=1}^p
E_X(X_{1j}X_{2j}X_{1l}X_{2l}X_{1r}X_{2r}X_{1s}X_{2s})\\
&=&9p^{-2}\sum_{j=1}^p\sum_{l=1}^p\sum_{r=1}^p \sum_{s=1}^p
H_{jlrs}^2,
\end{eqnarray*}
since $E_0(Y_1^4Y_2^4)=E_0^2(Y_1^4)=9,$ where we set
$$
H_{jlrs}\eq E_X(X_{1j}X_{1l}X_{1r}X_{1s})=\begin{cases} E_X(X_1^4),&
j=l=r=s,\\
1,& j=l\not=r=s\ \text{or}\ j=r\not=l=s\ \text{or}\
j=s\not=r=l,\\
0,& \text{otherwise} .
\end{cases}
$$
As a consequence, we get
$$
E_0(K^4(Z_1,Z_2))\le 9p^{-1}b_4^2+27\ ,
$$
where $b_4\eq \sup_{i}E(X_{i1}^4)$.
 By {\bf B1}, the second convergence in \nref{Unorm} holds
true. Thus, Theorem \ref{TU2} (i) follows.\\

Let us now evaluate the type II errors  under $P_\t$. Recall that
by \nref{Mod},
$$
Y_i=\xi_i+v_i,\quad v_i=\sum_{j=1}^p\t_jX_{ij},\quad \xi_i\sim
\CN(0,1)\quad iid.
$$
Observe that $E_\t Y_iX_{ij}=\t_j$ and set
$$
K_\t(Z_1,Z_2)=p^{-1/2}\sum_{j=1}^p(Y_1X_{1j}-\t_j)(Y_2X_{2j}-\t_j).
$$
Consider the representation
$$
K(Z_1,Z_2)=K_\t(Z_1,Z_2)+\delta(Z_1)+\delta(Z_2)+h(\t)
$$
where
$$
\delta(Z_i)=p^{-1/2}\sum_{j=1}^p(Y_iX_{ij}-\t_j)\t_j\ ,\quad \
h(\t)=p^{-1/2}\sum_{j=1}^p\t_j^2\ .
$$
Observe that the kernel $K_\t(Z_1,Z_2)$ is symmetric and degenerate
under $P_\t$, i.e.,
$$
E_\t^{Z_1}K_\t(Z_1,Z_2)=E_\t^{Z_2}K_\t(Z_1,Z_2)=0.
$$
The terms  $K_\t(Z_1,Z_2)$, $\delta(Z_1)$, and $\delta(Z_2)$ are centered and
uncorrelated under $P_\t$. As a consequence, we derive that
\begin{eqnarray}
E_\t(K(Z_1,Z_2))&=&p^{-{1/2}}\|\theta\|^2,\label{esperance_K}\\
\Var_\t(K(Z_1,Z_2))&=&\Var_\t(K_\t(Z_1,Z_2))+\Var_\t(\delta(Z_1))+
\Var_\t(\delta(Z_2))\ .\label{majoration_variance_globale_t1}
\end{eqnarray}
Let us compute the variances. Let $\delta_{ij}$ be the Kronecker function. Using
the representation
\begin{eqnarray*}
K_\t(Z_1,Z_2)&=&p^{-1/2}\sum_{j=1}^p
\l(\xi_1X_{1j}+\sum_{r=1}^p\t_r(X_{1r}X_{1j}-\delta_{rj})\r)\\
&&\times
\l(\xi_2X_{2j}+\sum_{s=1}^p\t_s(X_{2s}X_{2j}-\delta_{sj})\r),
\end{eqnarray*}
 we derive that
\begin{eqnarray*}
E_\t^X(K_\t(Z_1,Z_2))&=&p^{-1/2}\sum_{j=1}^p \sum_{r=1}^p
\sum_{s=1}^p\t_r\t_s(X_{1r}X_{1j}-\delta_{rj})(X_{2s}X_{2j}-\delta_{sj}),\\
\end{eqnarray*}
 Denoting
$H_{rsj}=(X_{1r}X_{1j}-\delta_{rj})(X_{2s}X_{2j}-\delta_{sj})$
observe that
$$
\Var_XE_\t^X(K_\t(Z_1,Z_2))=p^{-1}\sum_{j=1}^p \sum_{r=1}^p
\sum_{s=1}^p\sum_{l=1}^p \sum_{u=1}^p \sum_{v=1}^p\t_r\t_s\t_u\t_v
E_X(H_{rsj}H_{uvl}).
$$
Note that
$$
E_X(H_{rsj}H_{uvl})= D_{rujl}D_{svjl},
$$
where (we omit the first index $i=1,2$ in $X_{ij}$)
\begin{eqnarray*}
D_{rujl}&=&E_X\l((X_{r}X_{j}-\delta_{rj})(X_{u}X_{l}-\delta_{ul})\r),\\
D_{svjl}&=&E_X\l((X_{s}X_{j}-\delta_{sj})(X_{v}X_{l}-\delta_{vl})\r).
\end{eqnarray*}
Observe that
$$
D_{rujl}=\begin{cases}
1,&r=l\not=u=j\quad \text{or}\quad r=u\not=j=l,\\
b_4-1, &r=u=j=l,\\
0,& \text{otherwise}.
\end{cases}
$$
We obtain
\begin{eqnarray*}
&&\Var_XE_\t^X(K_\t(Z_1,Z_2))=p^{-1}\sum_{j=1}^p \sum_{r=1}^p
\sum_{s=1}^p\sum_{l=1}^p \sum_{u=1}^p \sum_{v=1}^p\t_r\t_s\t_u\t_v
D_{rujl}D_{svjl}\\
&&= \frac{(b_4-1)^2}{p}\sum_{r=1}^p\t_r^4+\frac{2b_4-1}{p}\sum_{1\le
r,s\le p,\ r\not=s}\t_r^2\t_s^2 +\frac{1}{p}\sum_{1\le j, r,s\le p,\
 j\not=r,
j\not=s}\t_r^2\t_s^2\\
&&=
O[\sum_{j=1}^p\t_j^4]+O[(\sum_{j=1}^p\t_j^2)^2]=O(\|\t\|^4).
\end{eqnarray*}
We now compute $E_X[\Var_\t^X(K_\t(Z_1,Z_2))]$.
\begin{eqnarray*}
\Var_\t^X(K_\t(Z_1,Z_2))&=&p^{-1}\sum_{1\le j,l\le p }
X_{1j}X_{2j}X_{1l}X_{2l}\\
&+& p^{-1}\sum_{j=1}^p\sum_{l=1}^p\sum_{r=1}^p\sum_{s=1}^p
X_{1j}X_{1l}\theta_r\theta_s(X_{2r}X_{2j}-\delta_{jr})(X_{2s}X_{2l}-\delta_{sl}
)\\
&+& p^{-1}\sum_{j=1}^p\sum_{l=1}^p\sum_{r=1}^p\sum_{s=1}^p
X_{2j}X_{2l}\theta_r\theta_s(X_{1r}X_{1j}-\delta_{jr})(X_{1s}X_{1l}-\delta_{sl}
)\ .
\end{eqnarray*}
Let us take the expectation with respect to $X$. By  Assumption {\bf A2}, we
have
\begin{eqnarray*}
E_X[\Var_\t^X(K_\t(Z_1,Z_2))]&=&1+
2p^{-1}\sum_{j=1}^p\sum_{r=1}^p\sum_{s=1}^p\theta_r\theta_sE_X[(X_{2r}X_{2j}
-\delta_ {jr})(X_{2s}X_{2j}-\delta_{sj})]\\
& \leq & 1+ 2\sum_{r=1}^p b_4 \theta_r^2 = 1+ O(\|\theta\|^2)= O(1+
\|\theta\|^4)
\end{eqnarray*}
Since
$$
\Var_\t(K_\t(Z_1,Z_2))=E_X\Var_\t^X(K_\t(Z_1,Z_2))+\Var_XE_\t^X(K_\t(Z_1,Z_2)),
$$
we get
\begin{eqnarray}
\Var_\t(K_\t(Z_1,Z_2))=O(1+ \|\t\|^4).\label{majoration_variance_t1}
\end{eqnarray}
Similarly for $i=1,2$, we compute the variance of $\delta(Z_i)$.
$$
\delta(Z_i)=p^{-1/2}\sum_{j=1}^p\t_j
\l(\xi_iX_{ij}+\sum_{l=1}^p\t_l(X_{ij}X_{il}-\delta_{jl})\r),
$$
and we have (we omit the index $i=1,2$)
\begin{eqnarray*}
E_\t^X(\delta(Z))&=&p^{-1/2}\sum_{j=1}^p\sum_{l=1}^p\t_j\t_l(X_{j}X_{l}-
\delta_{jl}),\\
\Var_\t^X(\delta(Z))&=&p^{-1}\sum_{j=1}^p\sum_{l=1}^p\t_j\t_{l}X_{j}X_{l},\quad
E_X\Var_\t^X(\delta(Z))=p^{-1}\|\t\|^2,\\
\Var_XE_\t^X(\delta(Z))&=&p^{-1}\sum_{j=1}^p\sum_{l=1}^p\sum_{r=1}^p\sum_{s=1}^p
\t_j\t_l\t_r\t_sD_{jrls}\ ,
\end{eqnarray*}
where $D_{jrls}$ was previously defined and upper bounded. This yields
\begin{eqnarray}\label{majoration_variance2t1}
\Var_XE_\t^X(\delta(Z))=p^{-1}\l((b_4-1)\sum_{j=1}^p\t_j^4+2\|\t\|^4\r)
=O(\|\t\|^4/p).
\end{eqnarray}

Combining \nref{esperance_K}, \nref{majoration_variance_globale_t1},
\nref{majoration_variance_t1}, and \nref{majoration_variance2t1} we obtain,
for
$r^2np^{-1/2}\to\infty$ and $p=o(n^2)$,
\begin{eqnarray*}
E_\t(t_1)&=&\sqrt{N}E_\t(K(Z_1,Z_2))=\sqrt{N}h(\t)\sim
n(2p)^{-1/2}\|\t\|^2\ge \frac{n}{\sqrt{2p}}r^2\to\infty,\\
\Var_\t(t_1)&=&\Var_\t(K_{\theta}(Z_1,Z_2))+
\frac{n^3}{N}\Var_\t(\delta(Z_1))=O(1+\|\t\|^4)+O(n\|\t\|^4/p)\\&=&
o\l(\l(E_\t(t_1)\r)^2\r).
\end{eqnarray*}
Applying Chebyshev's inequality as in the proof of Theorem \ref{TU1} allows to
conclude.
\endproof

\subsection{Higher Criticism Tests}\label{HCtest}

\subsubsection{Type I errors}\label{T1_HC}

The variables $X_1,...,X_p,Y$ are independent under $P_0$ and
$(X_j,a)/\|a\|\sim \CN(0,1)$ for any $a\in \R^p,\ a\not= 0$ under
{\bf A3}. Thus we have
\begin{eqnarray*}
P_0(y_1<t_1,...,y_p<t_p)&=&E_Y(P_0^Y((X_1,Y)/\|Y\|<t_1,...,(X_p,Y)/\|Y\|<t_p)\\
&=& E_Y(\Phi(t_1)....\Phi(t_p)) =\Phi(t_1)....\Phi(t_p)
=P_0(y_1<t_1)...P_0(y_p<t_p).
\end{eqnarray*}
It follows that $y_j=(X_j,Y)/\|Y\|\sim \CN(0,1)$ and $y_1,\dots,y_p$
are i.i.d. under $P_0$. As a consequence, the random variables $q_i$
are independent uniformly distributed on $(0,1)$  under $P_0$. We
denote by $F_p(t)$ the empirical distribution of $(q_i)_{1\leq i\leq
p}$:
$$F_p(t)=\frac{1}{p}\sum_{i=1}^{p}\1_{q_i\leq t}\ .$$
Then, the normalized uniform empirical process is defined by
$$W_p(t)=\sqrt{p}\frac{F_p(t)-t}{\sqrt{t(1-t)}}\ .$$
Arguing as in Donoho and Jin~\cite{dj04}, we observe that $t_{HC}=\sup_{t\le
1/2}W_p(t)$.
It is stated in \cite{shorack}, Chapter 16 that
$$ \frac{\sup_{0\leq t\le 1/2}W_p(t)}{\sqrt{2\log\log p}}\rightarrow_P 1\ ,
\quad p\rightarrow \infty\ .$$ This proves the result.\endproof

\subsubsection{Type II errors}\label{T2_HC}

We define $H_{np}=(1+a)\sqrt{2\log\log p}$. Consider some
$\beta\in(1/2,1)$ and assume that $k\log(p)/n\rightarrow 0$. It is
sufficient to prove that for any $\delta_0>0$ arbitrarily small the
radius
\begin{eqnarray}\label{definition_radius}
r_{np}=(\varphi(\beta)+\delta_0)\sqrt{k\log(p)/n}
\end{eqnarray}
satisfies
\begin{equation}\label{objectif_higher_criticism}
 \beta(\psi^{HC},\Theta_k(r_{np}))\rightarrow 0\ .
\end{equation}

For any $\theta\in\Theta_k$, we set $\|\theta\|_{\infty}\eq
\sup_i|\t_i|$. In order to prove the convergence
\nref{objectif_higher_criticism}, we consider a partition of
$\Theta_k(r_{np})$:
\begin{eqnarray*}
\tilde{\Theta}^{(1)}_k(r_{np}) &\eq& \Theta_k(r_{np})\cap\l\{\theta\in\Theta_k ,
\|\theta\|^2\geq \frac{4k\log(p)}{n}\r\}\\
 \tilde{\Theta}^{(2)}_k(r_{np})& \eq   & \Theta_k(r_{np})\cap
[\tilde{\Theta}^{(1)}_k(r_{np})]^c
\cap \l\{\theta\in\Theta_k\ , \|\theta\|^2_{\infty}\geq
\frac{4\log(p)}{n}\r\}\\
\tilde{\Theta}^{(3)}_k(r_{np})& \eq  & \Theta_k(r_{np})\cap
[\tilde{\Theta}^{(1)}_k(r_{np})]^c
\cap [\tilde{\Theta}^{(2)}_k(r_{np})]^c\ .
\end{eqnarray*}
The sets $\tilde{\Theta}^{(1)}_k(r_{np})$ and $\tilde{\Theta}^{(2)}_k(r_{np})$
contain the parameters $\theta$ whose $l_2$ or $l_{\infty}$ norms
are large, while the set $\tilde{\Theta}^{(3)}_k(r_{np})$ contains the remaining
parameters.

\begin{proposition}\label{prop_control_large}
Consider the set of parameters $\tilde\Theta^{(4)}_k$ defined by
$$
\tilde\Theta^{(4)}_k \eq \l\{\theta\in\Theta_k,\quad
\frac{\|\theta\|^2_{\infty}}{1+\|\theta\|^2}\geq \frac{3\log(p)}{n}\r\}\ . $$
Let us introduce the statistic $t_{\max}$  and the corresponding test
$\psi^{\max}$ defined by
\begin{eqnarray*}
 t_{\max}\eq  (pq_{(1)})^{-1/2}- (pq_{(1)})^{1/2}\leq
t_{HC}\ ,\  \quad  \psi^{\max}\eq \1_{t_{\max} > H_{np}}\ .
\end{eqnarray*}
We have $\beta(\psi^{\max},\tilde\Theta^{(4)}_k)\rightarrow 0$.
\end{proposition}
It follows that $\beta(\psi^{HC},\tilde\Theta^{(4)}_k)\rightarrow 0$.
Observe that $$\tilde{\Theta}^{(1)}_k(r_{np})\subset \l\{\theta\in\Theta_k\
,\quad
\frac{\|\theta\|^2}{1+\|\theta\|^2}\geq \frac{4k\log(p)/n}{1+4k\log(p)/n}\r\}\
.$$
Since $k\|\t\|_{\infty}^2\geq \|\t\|^2$ and since $k\log(p)/n$ converges to
$0$, it follows that $\tilde{\Theta}^{(1)}_k(r_{np})\subset
\tilde\Theta^{(4)}_k$ for $n$ large enough. Thus, we get
$\beta(\psi^{HC},\tilde\Theta^{(1)}_k(r_{np}))\rightarrow 0$.\\

Let us turn to $\tilde{\Theta}^{(2)}_k(r_{np})$. For any $\t\in
\tilde{\Theta}^{(2)}_k(r_{np})$, we have
$$\frac{\|\t\|^2_{\infty}}{1+\|\t\|^2}\geq \frac{4\log(p)/n}{1+ 4k\log(p)/n}\ .
$$ This quantity is larger than $3\log(p)/n$ for $n$ large enough. We get
$\beta(\psi^{HC},\tilde\Theta^{(2)}_k(r_{np}))\rightarrow 0$.

\begin{proposition}\label{prop_control_restant}
Let us set $T_p= \sqrt{\log(p)}$ and $u>0$ such that
\begin{equation}\label{choice1}
u  =\begin{cases}2 \varphi(\b), &\b\in (1/2,3/4]\ ,\\
\sqrt{2}, & \b\in (3/4, 1)\ .
\end{cases}
\end{equation}
We consider the statistic $L(u)$ and the corresponding test $\psi^{L}$ defined
by
\begin{eqnarray*}
L(u)\eq \sum_{j=1}^p\frac{\1_{|y_{j}|>uT_p}-
2\Phi(-uT_p)}{\sqrt{2p\Phi(-uT_p})}\ , \quad
\quad \psi^{L}=\1_{L(u)\geq  H_{np}}\ .
\end{eqnarray*}
Then, $\beta(\psi^{L},\tilde\Theta^{(3)}_k(r_{np}))\rightarrow 0$. Moreover, we
have $L(u)\leq t_{HC}$, for $p$ large enough.
\end{proposition}
It follows from Proposition \ref{prop_control_restant} that
$\beta(\psi^{HC},\tilde\Theta^{(3)}_k(r_{np}))\rightarrow 0$
converges to $0$, which completes the proof.
\endproof

\subsubsection{Proof of Proposition \ref{prop_control_large}}
It follows directly from the definition \nref{definition_HC} that $t_{\max}\leq
t_{HC}$.
Consider the test $\psi^{'\max}$ defined by
\begin{eqnarray}\label{definition_psi'}
\psi^{'\max}=\1_{\|y\|_{\infty}\geq \sqrt{2.5\log(p)}}
\end{eqnarray}
If $\psi^{'\max}=1$, it follows that $q_{(1)}\leq
2\Phi(-\sqrt{2.5\log(p)})\leq 2p^{-5/4}$. Hence, we have
$t_{\max}\geq p^{1/8}/\sqrt{2}-\sqrt{2}p^{-1/8}$. For $p$ large
enough, this implies that $\psi^{\max}=1$. Consequently, we only
have to prove that $\beta(\psi^{'\max},\tilde\Theta^{(4)}_k)\rightarrow 0$.\\

Consider $\t\in \tilde\Theta^{(4)}_k$. By symmetry, we may assume that
$\|\t\|_{\infty}=|\theta_1|$. We use the following decomposition
$$\|Y\|y_1= \theta_1\|X_1\|^2+ (Y-\theta_1X_1,X_1)\ .$$
The random variables $\|Y\|^2/(1+\|\t\|^2)$ and $\|X_1\|^2$ have a
$\chi^2$ distribution with $n$ degrees of freedom. Since
$Y-\theta_1X_1$ is independent of $X_1$, the random variable
$(Y-\theta_1X_1,X_1/\|X_1\|)$ is normal with mean $0$ and variance
~\\$1+\sum_{i\neq 1}\t_i^2$.

With probability larger than $1-O(n^{-1}\vee \log^{-1}(p))$,  we
obtain
\begin{eqnarray*}
\|Y\|^2/n&\leq& (1+\|\theta\|^2)[1+o(n^{-1/4})] \\
 (1- o(n^{-1/4}))\leq  \|X_1\|^2/n &\leq & (1+o(n^{-1/4}))\\
|(Y-\theta_1X_1,X_1)|/\|X_1\|&\leq & (1+\sum_{i\neq
1}\t_i^2)^{1/2}\sqrt{2\log(\log(p))}\ .
\end{eqnarray*}
Thus, we get
$$|y_1|\geq
\frac{\sqrt{n}|\theta_1|}{(1+\|\theta\|^2)^{1/2}}[1-o(n^{-1/4)}]-O(\sqrt{
\log \log(p)}) \ ,$$ with probability larger than $1-O(n^{-1}\vee
\log^{-1}(p))$. Since $\t\in\tilde{\Theta}^{(4)}_k$, we have
$n|\theta_1|^2/(1+\|\theta\|^2) \geq 3\log(p)$ and the test
$\psi'_{\max}$ rejects with probability going to one. It follows
that $\beta(\psi^{'\max},\tilde\Theta^{(4)}_k)\rightarrow 0$.
\endproof

\subsubsection{Proof of Proposition \ref{prop_control_restant}}

{\bf Connection between $t_{HC}$ and $L(u)$}. Set $\hat{s}_u \eq
\sum_{i=1}^p\1_{|y_j|>uT_p}$. Observe that $q_{(\hat s_u)}\leq
P(|\CN(0,1)|>uT_p)\leq 1/2$ for $p$ large enough. If follows that
$$L(u)= \frac{\sqrt{p}[\hat{s}_u/p- 2\Phi(-uT_p)]}{\sqrt{2\Phi(-uT_p)}}\leq
\frac{\sqrt{p}[\hat{s}_u/p- q_{(\hat{s}_u)}]}{\sqrt{q_{(\hat{s}_u)}}}\leq
t_{HC}\ .$$

\noindent {\bf Power of $\psi^{L}$}. Under $P_{\theta}$,
$\|Y\|^2/(1+\|\theta\|^2)$ has a $\chi^2$ distribution with $n$
degrees of freedom. For any $\t\in\tilde\Theta_k^{(3)}(r_{np})$, we
have $\|\theta\|^2\leq 4k\log(p)/n=o(1)$. As a consequence, we have
$|\|Y\|^2-n|\leq 4k\log(p)+4\sqrt{n\log(n)}=o(n)$ with probability
larger than $1-O(1/n)$ uniformly over all
$\t\in\tilde\Theta_k^{(3)}(r_{np})$. Consider the event
$\CZ_{np,1}=\{|\|Y\|^2-n|\le H_n\}$, where $H_n=
4k\log(p)+4\sqrt{n\log(n)}=o(n)$.
 It is sufficient to prove that
\begin{equation}\label{choice}
\sup_{\t\in\tilde\Theta^{(3)}_k(r_{np})}P_\t(\CZ_{np,1}\cap\{L(u)\le
H_{np}\})\to 0.
\end{equation}

Consider $\t\in \tilde\Theta^{(3)}_k(r_{np})$. We can assume that
$\t_{k+1}=...=\t_p=0$. Then $Y=\sum_{j=1}^k \t_jX_j+\xi$ does not
depend on $X_{k+1},...,X_p$. Arguing as for the type I error, we
derive that $y_{k+1},...,y_p$ are independent standard Gaussian
variables and do not depend on $(y_1,...,y_k)$. We can write $
L(u)=L_1(u)+L_2(u)$, where
\begin{eqnarray*}
L_1(u)&=&\frac{\sum_{j=1}^k\l(
\1_{\{|y_j|>uT_p\}}-2\Phi(-uT_p)\r)}{\sqrt{2p\Phi(-uT_p)}},\\
L_2(u)&=&\frac{\sum_{j=k+1}^p\l(
\1_{\{|y_j|>uT_p\}}-2\Phi(-uT_p)\r)}{\sqrt{2p\Phi(-uT_p)}}.
\end{eqnarray*}
We find
$$
 E_\t(L_2(u))=0,\quad \Var_\t(L_2(u))=
\frac{2p\Phi(-uT_p)(1-2\Phi(-uT_p))}{2p\Phi(-uT_p)}\le 1\ ,
$$
which yields,
\begin{equation}\label{L2}
P_{\theta}(|L_2(u)|> H_{np})\to 0\ .
\end{equation}

In order to study the term $L_1(u)$, we will find a statistic
$\tilde L_1(u)$ such that  $P_{\t}[\tilde L_1(u)<L_1(u)]=1+o(1)$ uniformly over
$\Theta^{(3)}_k(r_{np})$. For such a $\tilde L_1(u)$, we will have
\begin{equation}\label{L2aa}
P_\t[L(u)\le H_{np}]\leq P_\t[L_1(u)\le 2H_{np}]+o(1)\le P_\t[\tilde L_1(u)\le
2H_{np}]+o(1).
\end{equation}

\noindent
{\bf Construction of $\tilde L_1(u)$}. Observe that under $P_\t$,
\begin{eqnarray*}
y_j&=&\l(\hat y_j\|\xi\|+n\t_j+\Delta_j\r) /\|Y\|\ , \\
\Delta_j&=&\sum_{l\neq j}^k\t_l(X_j,X_l)+ \l(\|X_j\|^2-n\r)\t_j,
 \quad j=1,...,k,
\end{eqnarray*}
where
$$
\hat y_j=(X_j,\xi)/\|\xi\| .
$$
We only need to consider $Z\in \CZ_{np,2}= \{\|\xi\|^2-n|<n^{2/3}\}$
since $P_\t(\CZ_{np,2})\to 1$. Set
$\CZ_{np,3}=\CZ_{np,1}\cap\CZ_{np,2}$. Thus, for
$\delta=\delta_{np}\rightarrow 0$
one has
\begin{eqnarray}
\{|y_j|>uT_p\}\cap \CZ_{np,3}&\supset& \{|n^{-1/2}\hat
y_j\|\xi\|+n^{-1/2}\Delta_j+n^{1/2}\t_j)|>uT_p(1+\delta)\}\cap
\CZ_{np,3}\nonumber\\
&\supset& \{sgn(\t_j)\hat
y_j(1-\delta)>uT_p(1+\delta)-n^{1/2}|\t_j|+|\tilde S_j|)\}\cap
\CZ_{np,3}\ \nonumber,\\\label{relation_inclusion}
\end{eqnarray}
where $\tilde S_j=n^{-1/2}\Delta_j$.

\begin{lemma}\label{lemma_controle_deviation_S}
For any $T>0$ going to infinity and such that $T=o(\sqrt{n})$, we have
$$
\log(P_X(|\tilde S_j|>T\|\theta\|))\le -\frac 14
T^2(1+o(1))\ , $$
 uniformly over $\in\tilde\Theta^{(3)}_k(r_{np})$.
\end{lemma}

Taking $T=\sqrt{4\log(p)}$, we obtain
$$
P_X(|\tilde
S_j|>T\|\theta\|)=o(p^{-1}).
$$
We recall that  $\|\theta\|^2\leq 4k\log(p)/n=o(1)$ since
$\theta\in\tilde\Theta^{(3)}_k(r_{np})$. Hence, we get
$$
P_X\l[\max_{1\le j\le k}|\tilde S_j|>o(\sqrt{\log(p)})\r]=o(1),\quad
\text{uniformly over }\t \in\tilde\Theta^{(3)}_k(r_{np}).
$$
Combining this bound with  (\ref{relation_inclusion}), we obtain
that there exists an event $\CZ_{np,4}$ of probability tending to
one and a sequence $\delta=\delta_{np}\rightarrow 0$ such that
\begin{equation}
\{|y_j|>uT_p\}\cap \CZ_{np,4}\supset
 \{sgn(\t_j)\hat
y_j>uT_p(1+\delta)-(1-\delta)n^{1/2}|\t_j|\}\cap
\CZ_{np,4}\ . \label{relation_inclusion2}
\end{equation}
Observe that the random variables $\hat
y_j$ are independent standard normal.\\

Setting $\tilde u=u(1+\delta),\ \tilde
\rho_j=(1-\delta)n^{1/2}|\t_j|$ we define
$$
\tilde L_1(u)=\frac{\l(\sum_{j=1}^k\1_{\hat y_j>\tilde
uT_p-\tilde\rho_j}-2\Phi(-uT_p)\r)} {\sqrt{2p\Phi(-uT_p)}}.
$$
By (\ref{relation_inclusion2}), $\tilde L_1(u)$ satisfies
$P_\t[\tilde{L}_1(u)\leq L_1(u)]=1-o(1)$ uniformly over
$\tilde\Theta^{(3)}_k(r_{np})$. In view of \nref{L2aa}, in order to
complete the proof it suffices to show that
\begin{equation}\label{*1}
P_\t[\tilde L_1(u)\leq 2H_{np}]=o(1) \quad \text{uniformly over
$\tilde\Theta^{(3)}_k(r_{np})$.}
\end{equation}

\noindent
{\bf Control of $P_\t[\tilde L_1(u)\leq 2H_{np}]$.}
In order to evaluate this probability, recall that $\hat
y_j\sim\CN(0,1)$  i.i.d. under $P_\t$. Thus,
\begin{eqnarray*}
E_\t(\tilde L_1(u))&=&\frac{\l(\sum_{j=1}^k\Phi(-\tilde
uT_p+\tilde\rho_j)-2\Phi(-uT_p)\r)} {\sqrt{2p\Phi(-uT_p)}}\ ,\\
\Var_\t(\tilde L_1(u))&\le&\frac{\sum_{j=1}^k\Phi(-\tilde
uT_p+\tilde\rho_j)}{2p\Phi(-uT_p)}\ .
\end{eqnarray*}
By  Chebyshev's inequality, we get
\begin{eqnarray*}
P_\t(\tilde L_1(u)\le 2H_{np})&=& P_\t(E_\t(\tilde L_1(u))-\tilde
L_1(u)\ge E_\t(\tilde L_1(u))-2H_{np})\\
&\le& \frac{\Var_\t(\tilde L_1(u))}{(E_\t(\tilde
L_1(u))-2H_{np})^2}.
\end{eqnarray*}

\begin{lemma}\label{lemma_HC}
 There exists
$\eta>0$ such that, for $n,p$ large enough,
\begin{equation}\label{EL2}
\inf_{\t\in\tilde\Theta^{(3)}_k(r_{np})} \frac{\sum_{j=1}^k\Phi(-\tilde
uT_p+\tilde\rho_j)}{\sqrt{2p\Phi(-uT_p)}}\sim
\inf_{\t\in\tilde\Theta^{(3)}_k(r_{np})}
E_\t(\tilde L_1(u))>p^\eta.
\end{equation}

\end{lemma}

In the sequel, we denote by $A_p$ a log-sequence, i.e., a sequence
such that $A_p=(\log(p))^{c_p},\ |c_p|=O(1)$ as $p\to\infty$. Since
$u\in [0,\sqrt{2}]$, we have $p\Phi(-uT_p)\ge A_p$. Combining this
bound with Lemma \ref{lemma_HC} yields
$$
\Var_\t(\tilde L_1(u))=O\l(A_p E_\t(\tilde L_1(u))\r).
$$
Since $H_{np}=o(p^{\eta})$, this implies \nref{*1} and then  \nref{choice}.
\endproof

\subsubsection{Proof of Lemma \ref{lemma_controle_deviation_S}}

Let us bound the deviations of $\tilde
S_j$ by computing the
exponential moments of $\Delta_j$. For any $h$ such that $h^2\|\theta\|^2\leq
1/4$, we  have
\begin{eqnarray*}
E_X(\exp(h\Delta_j))&=&E_{X_j}E^{X_j}_X(\exp(h\Delta_j))\\
&=& E_{X_j}\l(\exp\l(h\t_j(\|X_j\|^2-n\r)
E^{X_j}_X \exp\l(h\sum_{l\not=j}\t_l(X_j,X_l)\r)\r)\\
&=&E_{X_j}\l(\exp\l(h\t_j(\|X_j\|^2-n)+\frac{h^2}{2}\|X_j\|^2
\sum_{l\not=j}\t_l^2\r)\r)\\
&=& \exp\l(-nh\t_j- \frac n2\log\l(1-
2h\t_j-h^2\sum_{l\not=j}\t_l^2\r)\r),
\end{eqnarray*}
as $2h\t_j+h^2\sum_{l\not=j}\t_l^2<1$.
Taking $h$ such that $h^2k\log(p)=o(n)$ enforces $h^2\|\theta\|^2=o(1)$.
Using the Taylor expansion of
the logarithm
$$
-hx-\frac 12\log(1- 2hx-h^2y^2)=\frac{1}{2}h^2(2x^2+y^2)(1+o(1)),\
h^2(2x^2+y^2)=o(1),
$$
we get
\begin{eqnarray}
E_X(\exp(h\Delta_j))&=&E_X(\exp(\sqrt{n}h\tilde
S_j))=\exp\l(\frac{n}{2}h^2(2\theta_j^2+
\sum_{l\not=j}\t_l^2)(1+o(1))\r) \nonumber\\
&\leq &\exp\l[nh^2\|\theta\|^2(1+o(1))\r]\label{exp_mom}
\end{eqnarray}
as $h^2\|\t\|^2=o(1)$. Take some $T>0$. Applying a standard technique based on
Markov's inequality yields
\begin{eqnarray*}
P_X(\tilde S_j>T\|\theta\|)&\le& E_X\l[\exp\l(\frac{T\tilde
S_j}{2\|\theta\|}-\frac{T^2}{2}\r)\r],\\
 P_X(-\tilde
S_j>T\|\theta\|)&\le& E_X\l[\exp\l(-\frac{T\tilde
S_j}{2\|\theta\|}+\frac{T^2}{2}\r)\r]\ .
\end{eqnarray*}
We get from \nref{exp_mom} that
$$
\log(P_X(|\tilde S_j|>T\|\theta\|))\le -\frac 14
T^2(1+o(1))\quad\text{if}\quad T^2=o(n)\text{ and } T\rightarrow \infty.
$$
\endproof

\subsubsection{Proof of Lemma \ref{lemma_HC}}

Recall that we consider  $r_{np}=(\varphi(\b)+
\delta_0)\sqrt{k\log(p)/n}$  with arbitrarily small $\delta_0>0$
(see \nref{definition_radius}). Recalling that $T_p=\sqrt{\log(p)}$,
we apply the results of Section \ref{extreme} for
$\delta=\delta_{np}>0,\ \delta_{np}=o(1)$, and
$$
T=\tilde uT_p,\quad \tilde u=(1+\delta)u,\quad
v=(1-\delta)(\varphi(\b)+\delta_0)<\tilde u,\quad t_0= v T_p,\quad
R=2T_p,
$$
since for $t_j=\tilde \rho_j$ one has
$$
\sum_{j=1}^kt_j^2=(1-\delta)^2n\sum_{j=1}^k\t_j^2\ge
(1-\delta)^2nr_{np}^2=kt_0^2.
$$

By the choice of $u,$ and  $v$, the relations \nref{L_ass} hold
true for $p$ large enough (see Remark \ref{R_ass}). Applying Lemmas
\ref{L_extr} and \ref{L_extr1}, we get
$$
\inf_{\t\in\tilde\Theta_k(r_{np})}\sum_{j=1}^k\Phi(-\tilde
uT_p+\tilde \rho_j)=k\Phi(-\tilde uT_p+t_0).
$$
We recall that $A_p$ denotes any $\log$-sequence.
Since $\Phi(-tT_p)=A_pp^{-t^2/2}$ for $t>0$, we have
\begin{eqnarray*}\nonumber
\inf_{\t\in\tilde\Theta_k(r_{np})}E_\t(\tilde
L_1(u))&=&\frac{k\l(\Phi(-\tilde
uT_p+t_0)-2\Phi(-uT_p)\r)}{\sqrt{2p\Phi(-uT_p)}}\sim
\frac{k\Phi(-\tilde uT_p+t_0)}{\sqrt{2p\Phi(-uT_p)}}\\
&=&\frac{k\Phi(-(\tilde u-v)T_p)}{\sqrt{2p\Phi(-uT_p)}}=A_p
p^{1/2-\b-(\tilde u-v)_+^2/2+u^2/4}.\nonumber
\end{eqnarray*}
In order to obtain \nref{EL2}, we have to check that there
exists $\eta>0$ such that, for $n, p$ large enough,
$$
G\eq \frac 12-\b-\frac{(\tilde u-v)_+^2}{2}+\frac{u^2}{4}\ge \eta.
$$

Let $\b\in (1/2,3/4]$. Recalling that $\varphi^2(\b)=2\b-1>0$ and
\nref{choice1} we see, that for $\delta=\delta_{np}=o(1)$ and
$\delta_0\in (0, \varphi(\b))$, one can find
$\eta=\eta(\b,\delta_0)>0$ such that
\begin{eqnarray*}
G&=&
-\frac{\varphi^2(\b)}{2}-\frac{(\varphi(\b)-\delta_0)^2}{2}+\varphi^2(\b)+o(1)\\
&=& \varphi(\b)\delta_0- \frac{\delta_0^2}{2}+o(1)\ge\eta+o(1).
\end{eqnarray*}

Let us now consider $\b\in (3/4,1]$. Recalling that
$\varphi(\b)=\sqrt{2}(1-\sqrt{1-\b})$ and \nref{choice1}, we see
that for $\delta=\delta_{np}=0(1)$ and $\delta_0\in (0, \sqrt{2-2\beta})$, one
can find
$\eta=\eta(\b,\delta_0)>0$ such that
\begin{eqnarray*}
G&=&\frac
12-\b-\frac{\l(\sqrt{2}-\sqrt{2}\l(1-\sqrt{1-\b}\r)-\delta_0\r)^2}{2}+\frac
12+o(1)\\
&=&1-\b-\l(\sqrt{1-\b}-\delta_0/\sqrt{2}\r)^2+o(1)\\
&=&\sqrt{2-2\b}\,\delta_0 -\frac{\delta_0^2}{2}+o(1)\ge\eta+o(1).
\end{eqnarray*}
The relation \nref{EL2} follows. \endproof

\subsection{Proof of Proposition \ref{prte_upper_bound_unknown}}

Under $H_0$, the distributions of the variables $(y_i)_{i=1,\ldots,p}$ do not
depend
on $\sigma^2$. As a consequence, $E_{0,\sigma}(\psi^{HC})=
E_{0,1}$. This last quantity has been shown to converge to $0$ in
Theorem \ref{TU4}. Hence, we get $\alpha^{un}(\psi^{HC})=o(1)$.\\

Let us turn to the type II error probability. We consider the model
$Y_i=\sum_{j=1}^p\theta_{j}X_{ij}+\xi_i$ where $\Var(\xi_i)=\sigma^2$. Dividing
this equation by $\sigma$, we obtain the model:
$$Y'_i=\sum_{j=1}^p(\theta_{j}/\sigma) X_{ij}+\xi'_i\ ,$$
where $\Var(\xi'_i)=1$. The statistic $t_{HC}$ is exactly the same
for the data $Z=(Y,X)$ and $Z'=(Y',X)$. Consequently, we obtain
$E_{\theta\sigma,\sigma}(1-\psi^{HC})= E_{\theta,1}(1-\psi^{HC})$.
It remains to use the bound on $E_{\theta,1}(1-\psi^{HC})$ from
Theorem \ref{TU4}.\endproof

\section{Appendix: Technical results}\label{SS1}

\subsection{Thresholds}\label{Tr} Take the thresholds $T=T_{j}$
satisfying
$$
T_j=\frac{a_j}{2}+\frac{\log(h^{-1})}{a_j}.
$$
Define
$a_j= x_j\sqrt{\log(p)}$, $\tau_j=T_j/\sqrt{\log(p)}$, and $h=
p^{-\b}$.
Then, we have $ \tau_j=x_j/2+\b/x_j$.\\

 If for some $\delta_0>0$,
$
x_j+\delta_0<\varphi_2(\b)\eq\sqrt{2}(1-\sqrt{1-\b})\leq \varphi(\beta)\ ,$
then there exists $\delta_1>0$ such that $\tau_j>\sqrt{2}+\delta_1
$. For such a $x_j$, we derive that
\begin{equation}\label{Phi0}
pT_j^r\Phi(-T_j)=o(1),\quad \forall\ r>0\ .
\end{equation}
In particular, if $x_j=o(1)$, then $\tau_j\to\infty$ and (\ref{Phi0}) holds.\\

\noindent
For any $\delta>0$, we have
\begin{equation}\label{Phi1}
\Phi(-T_j)\asymp  h\Phi(-T_j+a_j)\quad \text{for} \quad
\tau_j>x_j+\delta\ .
\end{equation}
This holds if $x_j<\varphi_2(\b)\le \sqrt{2}$.

\subsection{Norms $\|X_j\|$ and scalar products
$(X_j,X_l)$}\label{scalar}

Clearly,
$$
E(\|X_j\|^2)=n,\quad E(X_j,X_l)=0\ ,\quad \Var(X_j,X_l)=n.
$$
By Assumption ${\bf B1}$, there exists $D>0$ such that $\sup_{j\neq
l}\Var(X_j,X_l)\leq nD$ and $\sup_{j}\Var(\|X_j\|^2)\leq nD$.

\begin{lemma}\label{large_div} Let $U_j$ be a random
variable distributed as $X_{ij}$.

(1) Assume that there exists $h_0>0$ such that $\sup_{1\leq j\leq
l\leq p}E(e^{hU_jU_l})<\infty$ for any $|h|<h_0$. Then, for any
sequence $t=t_n$ such that $t=o(\sqrt{n})$ and $t\sqrt{n}\to\infty$,
\begin{eqnarray*}
P(|\|X_j\|^2-n|>t\sqrt{n})&\le& \exp[-t^2/(2D)(1+o(1))],
\end{eqnarray*}
and
\begin{eqnarray*}
P(|(X_j,X_l)||>t\sqrt{n})&\le& \exp[-t^2/(2D)(1+o(1))],
\end{eqnarray*}

(2) Assume that $E(|X|^m)<\infty$, for some $m>2$. Then there exists
$C_m<\infty$ such that
$$
P(|\|X_j\|^2-n|>t\sqrt{n})\le C_mt^{-m/2},\quad
P(|(X_j,X_l)|>t\sqrt{n})\le C_mt^{-m}.
$$

\end{lemma}
{\bf Proof} follows from the standard arguments based on the moment
inequalities and exponential inequalities. If $EZ=0,\
\Var(Z)=1,\ E(e^{h_0Z})<\infty$, then $\log(Ee^{hZ})=h^2/2(1+o(1))$
as $h\to 0$. Hence, we take $h=t/\sqrt{n}=o(1)$ for the study of the
exponential moments of $S_n=\sum_{i=1}^nZ_i$.

\begin{corollary}\label{C1} $\ $

\noindent(1) Let $\log(p)=o(n)$ and the assumptions Lemma
\ref{large_div} (1) hold true. Then, for any $B>2$, one has
\begin{eqnarray*}
P_X(\max_{1\le j\le p}|\|X_j\|^2-n|>\sqrt{BDn\log(p)}\,)=o(1),\\
P_X(\max_{1\le j<l\le p}|(X_j,X_l)|>\sqrt{2BDn\log(p)}\,)=o(1).
\end{eqnarray*}
(2) Let $p=o(n^{m/4})$ and the assumptions Lemma \ref{large_div} (2)
hold true. Then, for any sequence $v_n$ going to infinity, one has
\begin{eqnarray*}
P_X(\max_{1\le j\le p}|\|X_j\|^2-n|>\sqrt{n}p^{2/m}v_n)=o(1),\\
P_X(\max_{1\le j<l\le p}|(X_j,X_l)|>\sqrt{n}p^{2/m}v_n)=o(1).
\end{eqnarray*}
(3) Under assumptions (1) or (2) uniformly in $1\le j<l\le p$ in
$P_X$-probability, one has $a_j\sim a=b\sqrt{n},\quad x_j\sim x$,
i.e., for any $\delta>0$,
$$
P_X(\max_{1\le j\le p}|(a_j/b\sqrt{n})-1|>\delta)\to 0,\quad
P_X(\max_{1\le j\le p}|(x_j/x)-1|>\delta)\to 0.
$$
\end{corollary}

\subsection{Expansion of $\Phi(t)$}\label{exp_phi} Let $\Phi(t)$ be
the standard Gaussian cdf and $\phi(t)$ be the standard Gaussian
pdf.

\begin{lemma}\label{L4} Let $\delta\to 0,\
t\delta=O(1)$. Then
$$
\Phi(t+\delta)=\Phi(t)+\delta\phi(t)+O\l(\delta^2(|t|+1)\phi(t)\r).
$$
\end{lemma}
{\bf Proof} follows from the Taylor expansion and the properties of
$\phi(t)$.
\endproof

Observe that for any $b\in \R$ there exists $C=C(b)>0$ such that
$(|t|+1)\phi(-t)\le C(b)\Phi(-t)$ as $t\le b$. It follows
from Lemma \ref{L4} that as $\delta\to 0,\ t\delta=O(1),\ t\le B$
for some $B\in\R$, then
$$
\Phi(-t+\delta)=\Phi(-t)(1+O(\delta^2))+\delta\phi(t).
$$

\subsection{Tails  of correlated vectors}\label{prob_corr}

\begin{lemma}\label{L3} Let $(X,Y)$ be the Gaussian random two-dimensional
vector,
$$EX=EY=0,\ \Var(X)=\Var(Y)=1,\ \Cov(X,Y)=r.$$ Let $t_1\asymp
t_2\to\infty, \ rt_1=o(1)$. Then
$$
P(X>t_1,Y>t_2)=
\Phi(-t_1)\Phi(-t_2)\l(1+O(r^2)\r)+r\phi(t_1)\phi(t_2).
$$
\end{lemma}
{\bf Proof}. Observe that the conditional distribution $\CL(Y|X=x)$
is Gaussian $\CN(m(x),\sigma^2(x))$ with $m(x)=rx,\
\sigma^2(x)=1-r^2$. Therefore
$$
P(X>t_1,Y>t_2)=\int_{t_1}^\infty P(Y>t_2|X=x)d\Phi(x)=
\int_{t_1}^\infty\Phi\l(\frac{-t_2+rx}{\sqrt{1-r^2}}\r)d\Phi(x).
$$
Setting $h=|r|^{-1}$, observe that
$$
\int_{h}^\infty\Phi\l(\frac{-t_2+rx}{\sqrt{1-r^2}}\r)d\Phi(x)\le
\Phi(-h)=o\l(  r^2 \Phi(-t_1)\Phi(-t_2)\r) .
$$
It is sufficient to study the integral over the interval
$\Delta=[t_1,h]$. For $x\in \Delta$, we have
$$
\frac{-t_2+rx}{\sqrt{1-r^2}}=-t_2+\delta(x),\quad \delta(x)=
rx+O(r^2t_2+|r^3x|))=O(1).
$$
Applying Lemma \ref{L4},   we have
\begin{eqnarray*}
\lefteqn{\int_{\Delta}\Phi\l(\frac{-t_2+rx}{\sqrt{1-r^2}}\r)d\Phi(x)}&&\\
&=&
\Phi(-t_2)\l(\Phi(-t_1)-\Phi(-h)\r)\l(1+O(r^2)\r)
+r\phi(t_2)\int_\Delta xd\Phi(x)\\
 &=&
\Phi(-t_1)\Phi(-t_2)\l(1+O(r^2)\r)+r\phi(t_1)\phi(t_2),
\end{eqnarray*}
since $\int_\Delta
xd\Phi(x)=\phi(t_1)-\phi(h)=\phi(t_1)+o(r^2\Phi(-t_1))$. \endproof

\subsection{A minimization problem}\label{extreme}

Let $f(t)$ be a function defined on the interval $t\in [0,R]$.
Consider the minimization problem
\begin{equation}\label{extr}
F_k(t_0)=\inf\sum_{j=1}^k f(t_j)\quad\text{subject to}\quad
\sum_{j=1}^k t_j^2\ge kt_0^2,\quad t_j\in [0,R].
\end{equation}

\begin{lemma}\label{L_extr} Assume that there exists $\la>0$ such that
$$
\inf_{t\in [0,R]}(f(t)-\la t^2)=f(t_0)-\la t_0^2.
$$
Then $ F_k(t_0)=kf(t_0). $
\end{lemma}
{\bf Proof}. We have, for any $(t_1,...,t_k)$ such that $t_j\in
[0,R]$, $\sum_{j=1}^k t_j^2\ge kt_0^2$,
\begin{eqnarray*}
\sum_{j=1}^kf(t_j)&\ge& \sum_{j=1}^kf(t_j)-\la \l(\sum_{j=1}^k
t_j^2-kt_0^2\r)=\sum_{j=1}^k\l(f(t_j)-\la t_j^2\r)+k\la t_0^2\\
&\ge& k\l(f(t_0)-\la t_0^2\r)+k\la t_0^2= k f(t_0).\quad
\endproof
\end{eqnarray*}

We apply Lemma \ref{L_extr} to the function $f(x)=\Phi(-T+x)$. Let
$\phi(x)=\Phi^{\prime}(x)$ stand for the standard Gaussian pdf.

\begin{lemma}\label{L_extr1} Let $f(t)=\Phi(-T+t)$. Suppose
\begin{equation}\label{L_ass}
t_0>0,\quad T>t_0+\frac{2}{t_0},\quad
 T<R\le \l(\frac{t_0}{\phi(-T+t_0)}\r)^{1/2}.
\end{equation}
Take $\la=\phi(-T+t_0)/2t_0$. Then the assumptions of Lemma
\ref{L_extr} are fulfilled, i.e.,
$$
\inf_{0\le t\le R}(f(t)-\la t^2)=f(t_0)-\la t^2_0.
$$
\end{lemma}
{\bf Proof}. Denote $g(t)=\Phi(-T+t)-\la t^2$. By the choice of $\la$
we have $g^{\prime}(t_0)=0$. Let us consider the second derivative,
$$
g^{\prime\prime}(t)=(T-t)\phi(-T+t)-2\la=(T-t)\phi(-T+t)-\phi(-T+t_0)/t_0.
$$
Observe that the function $-x\phi(x)$ is positive for $x<0$,
increases for $x\in (-\infty,-1)$ and decreases for $x\in (-1,1)$;
$\lim_{x\to-\infty}\phi(x)=0=0\phi(0)$,
$$
g^{\prime\prime}(t_0)=(T-t_0-t_0^{-1})\phi(-T+t_0)>0.
$$
Consequently, there exist two points $t_1, t_2$ such that $t_1<t_0<t_2<T$,
$$
g^{\prime\prime}(t_1)=g^{\prime\prime}(t_2)=0,\quad
g^{\prime\prime}(t)<0\quad \text{as} \quad t<t_1\quad\text{and}\quad
t>t_2.
$$
The function $g(t)$ is therefore convex on $[t_1,t_2]$, concave on
$(-\infty,t_1]$ and on $[t_2,\infty)$, and $t_0$ is the point of a
local minimum of $g(t)$. By the concavity, this yields that the
global minimum of $g(t)$ at $t\in [0,R]$ is achieved either at
$t=t_0$ or at the ends of the interval $[0,R]$. Therefore we only
need to show that $g(0)>g(t_0)$ and $g(R)>g(t_0)$.

In order to verify the first inequality, observe that $g(0)>0$.
Recalling the well known inequality:
$$
 \Phi(-y)< \frac{1}{y}\phi(-y),\ \forall y>0\ ,
$$
we get
$$
g(t_0)=\Phi(-T+t_0)-t_0\phi(-T+t_0)/2<\phi(-T+t_0)
\l(\frac{1}{T-t_0)}-\frac{t_0}{2}\r)<0\ ,
$$
because $T>t_0+2t_0^{-1}$.

The second inequality follows from the relation
$$
g(R)=\Phi(-T+R)-\frac{R^2\phi(-T+t_0)}{2t_0}>
\frac{t_0-R^2\phi(-T+t_0)}{2t_0}>0,
$$
in view of the assumption on $R$. \endproof

\begin{remark}\label{R_ass} {\rm
Observe that if  $0<v<u<b$ and
$$
T=uT_p,\quad  t_0=vT_p,\quad  R=bT_p\ ,
$$
where $T_p$ is large enough, then assumptions \nref{L_ass} hold.}
\end{remark}

\end{document}